\documentclass[12pt]{article}

\usepackage{amssymb}

\newtheorem{thm}{Theorem}[section]
\newtheorem{defn}[thm]{Definition}
\newtheorem{prop}[thm]{Proposition}
\newtheorem{cor}[thm]{Corollary}
\newtheorem{lemma}[thm]{Lemma}
\newtheorem{rema}[thm]{Remark}

\newcommand{\halmos}{\rule{1ex}{1.4ex}}

\newcommand{\C}{\mathbb{C}}

\newcommand{\R}{\mathbb{R}}
\newcommand{\Z}{\mathbb{Z}}

\newcommand{\nn}{\nonumber \\}

 \newcommand{\pf}{{\it Proof.}\hspace{2ex}}
 
 \newcommand{\epfv}{\hspace*{\fill}\mbox{$\halmos$}\vspace{1em}}

\newcommand{\wt}{\mbox{\rm wt}\ }
\newcommand{\swt}{\mbox{\rm {\scriptsize wt}}\ }

\title{ {\bf Differential equations and  intertwining operators} }
\date{}
\author{Yi-Zhi Huang}

 \makeatletter
\newlength{\@pxlwd} \newlength{\@rulewd} \newlength{\@pxlht}
\catcode`.=\active \catcode`B=\active \catcode`:=\active
\catcode`|=\active
\def\sprite#1(#2,#3)[#4,#5]{
   \edef\@sprbox{\expandafter\@cdr\string#1\@nil @box}
   \expandafter\newsavebox\csname\@sprbox\endcsname
   \edef#1{\expandafter\usebox\csname\@sprbox\endcsname}
   \expandafter\setbox\csname\@sprbox\endcsname =\hbox\bgroup
   \vbox\bgroup
  \catcode`.=\active\catcode`B=\active\catcode`:=\active\catcode`|=\active
      \@pxlwd=#4 \divide\@pxlwd by #3 \@rulewd=\@pxlwd
      \@pxlht=#5 \divide\@pxlht by #2
      \def .{\hskip \@pxlwd \ignorespaces}
      \def B{\@ifnextchar B{\advance\@rulewd by \@pxlwd}{\vrule
         height \@pxlht width \@rulewd depth 0 pt \@rulewd=\@pxlwd}}
      \def :{\hbox\bgroup\vrule height \@pxlht width 0pt depth
0pt\ignorespaces}
      \def |{\vrule height \@pxlht width 0pt depth 0pt\egroup
         \prevdepth= -1000 pt}
   }
\def\endsprite{\egroup\egroup}
\catcode`.=12 \catcode`B=11 \catcode`:=12 \catcode`|=12\relax
\makeatother

\def\hboxtr{\FormOfHboxtr} 
\sprite{\FormOfHboxtr}(25,25)[0.5 em, 1.2 ex] 

:BBBBBBBBBBBBBBBBBBBBBBBBB |
:BB......................B |
:B.B.....................B |
:B..B....................B |
:B...B...................B |
:B....B..................B |
:B.....B.................B |
:B......B................B |
:B.......B...............B |
:B........B..............B |
:B.........B.............B |
:B..........B............B |
:B...........B...........B |
:B............B..........B |
:B.............B.........B |
:B..............B........B |
:B...............B.......B |
:B................B......B |
:B.................B.....B |
:B..................B....B |
:B...................B...B |
:B....................B..B |
:B.....................B.B |
:B......................BB |
:BBBBBBBBBBBBBBBBBBBBBBBBB |

\endsprite

\sprite{\FormOfShboxtr}(25,25)[0.3 em, 0.72 ex] 

:BBBBBBBBBBBBBBBBBBBBBBBBB |
:BB......................B |
:B.B.....................B |
:B..B....................B |
:B...B...................B |
:B....B..................B |
:B.....B.................B |
:B......B................B |
:B.......B...............B |
:B........B..............B |
:B.........B.............B |
:B..........B............B |
:B...........B...........B |
:B............B..........B |
:B.............B.........B |
:B..............B........B |
:B...............B.......B |
:B................B......B |
:B.................B.....B |
:B..................B....B |
:B...................B...B |
:B....................B..B |
:B.....................B.B |
:B......................BB |
:BBBBBBBBBBBBBBBBBBBBBBBBB |

\endsprite

\begin{document}

\bibliographystyle{alpha}
\maketitle

\begin{abstract}
We show that if every module $W$ for a vertex operator algebra
$V=\coprod_{n\in\mathbb{Z}}V_{(n)}$
satisfies the condition $\dim W/C_{1}(W)<\infty$, where
$C_{1}(W)$ is the subspace of $W$ spanned by
elements of the form $u_{-1}w$ for $u\in V_{+}=\coprod_{n>0}V_{(n)}$
and
$w\in W$, then matrix elements of
products and iterates of intertwining operators satisfy certain
systems of
differential equations.  Moreover, for prescribed singular points,
there exist such systems of differential equations such that the
prescribed singular points are regular. The finiteness of the fusion
rules is an immediate consequence of a result used to
establish the existence of such systems. Using these systems of
differential equations and some additional reductivity
conditions, we prove that products of intertwining operators
for $V$ satisfy the convergence and extension property needed in the
tensor
product theory for $V$-modules.
Consequently, when a vertex operator algebra $V$ satisfies all the
conditions
mentioned above, we obtain
a natural structure of vertex tensor category
(consequently
braided tensor category)  on the category
of $V$-modules and a natural structure of
intertwining operator algebra on
the direct sum of all (inequivalent) irreducible $V$-modules.
\end{abstract}

\renewcommand{\theequation}{\thesection.\arabic{equation}}
\renewcommand{\thethm}{\thesection.\arabic{thm}}
\setcounter{equation}{0}
\setcounter{thm}{0}
\setcounter{section}{-1}

\section{Introduction}

In the present paper, we show that  for a vertex operator algebra
satisfying certain finiteness and reductivity conditions, matrix
elements of products and iterates of intertwining operators satisfy
certain systems of differential equations of regular singular points.
Similar results are also obtained by Nagatomo and Tsuchiya in
\cite{NT}. Using these differential equations together with the
results
obtained by Lepowsky and the author in \cite{HL1}--\cite{HL4}
\cite{H1} and by the author in \cite{H1} \cite{H3}--\cite{H6},
we construct braided tensor categories and intertwining operator
algebras.

In the study of conformal field theories associated to affine Lie
algebras (the Wess-Zumino-Novikov-Witten models or WZNW models
\cite{W})
and to Virasoro algebras (the minimal models \cite{BPZ}), the
Knizhnik-Zamolodchikov equations \cite{KZ} and the
Belavin-Polyakov-Zamolodchikov equations \cite{BPZ}, respectively,
play
a fundamental role. Many important results for these and related
theories, including the constructions of braided tensor category
structures and the study of properties of correlation functions, are
obtained using these equations (see, for example, \cite{TK},
\cite{KL},
\cite{V}, \cite{H3}, \cite{HL5} and \cite{EFK}).

More generally, the tensor product theory for the category of modules
for a vertex operator algebra was developed by Lepowsky and the author
\cite{HL1}--\cite{HL4} \cite{H1} and the theory of intertwining
operator
algebras was developed by the author \cite{H1} \cite{H3}--\cite{H6}.
These structures are essentially equivalent to chiral genus-zero
weakly
conformal field theories \cite{S1} \cite{S2} \cite{H4} \cite{H5}.  In
the
construction of these structures {}from representations of vertex
operator
algebras, one of the most important steps is to prove the
associativity
of intertwining operators, or a weaker version which,
in physicists' terminology, is called the
(nonmeromorphic) operator product expansion of chiral vertex
operators.
It was proved in \cite{H1} that if a vertex operator algebra $V$ is
rational in the sense of \cite{HL1}, every finitely-generated
lower-truncated generalized $V$-module is a $V$-module and products of
intertwining operators for $V$ have a convergence and extension
property (see Definition 3.2 for the precise description of the
property),
then the associativity of intertwining operators holds. Consequently
the
category of $V$-modules has a natural structure of vertex tensor
category (and braided tensor category) and the direct sum of
all (inequivalent) irreducible $V$-modules has a natural structure of
intertwining operator algebra.

The results above reduce the construction of vertex tensor categories
and intertwining operator algebras (in particular, the proof of the
associativity of intertwining operators) to the proofs of the
rationality of vertex operator algebras (in the sense of \cite{HL1}),
the condition on
finitely-generated lower-truncated generalized $V$-modules and the
convergence and extension property. Note that this rationality and the
condition on finitely-generated lower-truncated generalized
$V$-modules
are both purely representation-theoretic properties. These and other
related representation-theoretic properties have been discussed in a
number of papers, including  \cite{Z},
\cite{FZ},  \cite{H1},  \cite{DLM1}--\cite{DLM3},
\cite{L}, \cite{KarL}, \cite{GN}, \cite{B}, \cite{ABD} and \cite{HKL}.

The convergence and extension property, on the other hand, is very
different {}from these purely representation-theoretic
properties. Since this property is for all intertwining operators, it
is impossible to prove it, even only the convergence part, by direct
estimates.  In the work \cite{TK}, Tsuchiya and Kanie used the
Knizhnik-Zamolodchikov equations to show the convergence of the
correlation functions.  For all the concrete examples (see \cite{H2},
\cite{H3}, \cite{HL5}, \cite{HM1} and \cite{HM2}), the convergence and
extension property was proved using the particular differential
equations of regular singular points associated to the examples,
including the Knizhnik-Zamolodchikov equations and the
Belavin-Polyakov-Zamolodchikov equations mentioned above.  Also, since
braided tensor categories and intertwining operator algebras give
representations of the braid groups, {}from the solution to the
Riemann-Hilbert problem, we know that there must be some differential
equations such that the monodromies of the differential equations give
these representations of the braid groups. In particular, it was
expected that there should be differential equations satisfied by the
products and iterates of intertwining operators.  In \cite{Ne}, by
using the fact that any finite-dimensional quotient space of a space
of functions with an action of the derivative operators gives a system
of differential equations, Neitzke observed that some cofiniteness
conditions (including the one discussed in this paper) are related to
systems of differential equations. Such differential equations can be
used to define correlation functions.  But the crucial problem of
whether the products or iterates of intertwining operators or the
formal series obtained from sewing three point correlation functions
satisfy differential equations was not addressed and was unsolved.

In the present paper, we solve this problem by showing that if every
module $W$ for a vertex operator algebra
$V=\coprod_{n\in\mathbb{Z}}V_{(n)}$ satisfies the condition $\dim
W/C_{1}(W)<\infty$ where $C_{1}(W)$ is the subspace of $W$ spanned by
elements of the form $u_{-1}w$ for $u\in V_{+}=\coprod_{n>0}V_{(n)}$
and $w\in W$, then matrix elements of products and iterates of
intertwining operators satisfy certain systems of differential
equations.  Moreover, for any prescribed singular point, there exist
such systems of differential equations such that this prescribed
singular point is regular. The finiteness of the fusion rules is an
immediate consequence of a result used to establish the existence of
such systems.  Together with some other reductivity and finiteness
conditions, including the complete reducibility of generalized modules
and the finiteness of the number of equivalence classes of modules,
these systems of differential equations also give us the convergence
and extension property.  Consequently, if all the conditions mentioned
above are satisfied, we obtain vertex tensor category structures
(consequently braided tensor categories) on the category of
$V$-modules and intertwining operator algebra structures on the direct
sum of all (inequivalent) irreducible $V$-modules.

The results of the present paper have many applications. For example,
in \cite{H7}, they have been used by the author in the construction of
the chiral genus-one correlation functions and in the proof of the
duality properties and modular invariance of these functions.  In
\cite{Y}, they have been used by Yamauchi in the study of module
categories of simple current extensions of vertex operator algebras.
They have also been also used to prove one of the Moore-Seiberg formulas 
in the author's proof of the Verlinde conjecture in \cite{H8}.

In \cite{Z1} and \cite{Z}, Zhu introduced the first cofiniteness
condition (called the $C_{2}$-cofiniteness condition in this paper) in
the theory of vertex operator algebras.  The condition $\dim
W/C_{1}(W)<\infty$, which we shall call the $C_{1}$-cofiniteness
condition (or we say that $W$ is $C_{1}$-cofinite), is a very minor
restriction and is very easy to verify for familiar examples.  Note
that vertex operator algebras as modules for themselves always satisfy
this condition.  (A different but closely-related condition, called
the $C_{1}$-finiteness condition, was discussed in \cite{L}.) This
condition was first introduced by Nahm in \cite{N} in which it was
called ``quasi-rationality.''  In \cite{AN}, Abe and Nagatomo showed
that if all modules for a vertex operator algebra satisfies the
$C_{1}$-cofiniteness condition (called $B_{1}$ finiteness condition by
these authors) and some additional conditions on the vertex operator
algebras are satisfied, then the spaces of conformal blocks on the
Riemann sphere are finite-dimensional.  This finite-dimensionality is
an immediate consequence of some results used to establish the
existence of the systems of differential equations in the present
paper; in particular, only the $C_{1}$-cofiniteness condition on
modules is needed.

In \cite{GN}, Gaberdiel and Neitzke found a useful spanning set for a
vertex operator algebra satisfying Zhu's $C_{2}$-cofiniteness
condition and in \cite{B}, Buhl found useful generalizations of this
spanning set for weak modules for such a vertex operator
algebra. Using these spanning sets, Abe, Buhl and Dong \cite{ABD}
proved that for a vertex operator algebra $V$ satisfying $V_{(n)}=0$
for $n<0$, $V_{(0)}=\mathbb{C}\mathbf{1}$ and the $C_{2}$-cofiniteness
condition, any irreducible weak $V$-module is also $C_{2}$-cofinite
and the complete reducibility of every $\mathbb{N}$-gradable weak
$V$-module implies that every weak $V$-module is a direct sum of
$V$-modules.  In Section 3, we also show that for such a vertex
operator algebra, the complete reducibility of every $V$-module
implies that every generalized $V$-module is a direct sum of
$V$-modules. These results allow us to replace the conditions needed
in the applications mentioned above by some stronger conditions which
have been discussed extensively in a number of papers.

Note that intertwining operator algebras have a geometric
formulation in terms of the sewing operation in the moduli space
of spheres with punctures and local coordinates (see \cite{H4} and
\cite{H5}). In fact, this formulation gives genus-zero modular
functors and genus-zero weakly conformal field theories in the
sense of Segal \cite{S1} \cite{S2}, and, in particular, it gives
the sewing property. It can be shown that the sewing property
together with the generalized rationality property of intertwining
operators (see \cite{H6}) implies the factorization property for
conformal blocks on genus-zero Riemann surfaces with punctures.  In
\cite{NT}, for a vertex operator algebra satisfying Zhu's
$C_{2}$-cofiniteness condition, the condition that Zhu's algebra
is semisimple and finite-dimensional and a condition stating that
certain induced weak $V$-modules are irreducible, Nagatomo and
Tsuchiya prove, among many other things, the factorization
property of conformal blocks on genus-zero surfaces. Their
factorization property is formulated and proved using the same
method as in the work \cite{TUY} by Tsuchiya, Ueno and Yamada on
WZNW models. In particular, they prove the convergence of
conformal blocks near the boundary points of the moduli space
using certain holonomic systems of regular singularities near the
singularities. The factorization theorem corresponds to the 
construction of the sewing operation for a modular functor in 
our formulation. On the other hand,
the sewing property for the corresponding 
weakly conformal field theories (in particular, 
the associativity of intertwing
operators and the operator product expansion of 
intertwining operators) are not discussed in \cite{NT}. 

Though many results in the present paper and in \cite{NT} are
equivalent, the present work and the work \cite{NT} develop into
different directions. The present work proves and studies in details
the properties of correlation functions (for example, associativity
and commutativity of intertwining operators) while the work \cite{NT}
proves and studies in details the algebro-geometric properties of
chiral genus-zero modular functors and conformal field theories (for
example, the construction and study of sheaves of conformal blocks
over the compactified moduli spaces).

We assume that the reader is familiar with the basic notions,
notations and results in the theory of vertex operator algebras as
presented in \cite{FLM} and \cite{FHL}. We also assume that the reader
is familiar with the theory of intertwining operator algebras as
developed and presented by the author in \cite{H1}, \cite{H3},
\cite{H4}, \cite{H5} and \cite{H6}, based on the tensor product theory
developed by Lepowsky and the author in \cite{HL1}--\cite{HL4} and
\cite{H1}.

The present paper is organized as follows: The existence of
systems of differential equations  and the existence of  systems
with  regular prescribed singular points are established in
Section 1 and Section 2, respectively. In Section 3, we prove the
finiteness of the fusion rules and the convergence and extension
property. Consequently we obtain the vertex tensor category (and
braided tensor category) structures and intertwining operator
algebra structures. We also discuss conditions which imply some
representation-theoretic conditions needed in the main
application. 

\paragraph{Acknowledgment}
I am grateful to A. Tsuchiya for a brief description of the
conditions and conclusions in his joint work with
K. Nagatomo and to K. Nagatomo for discussions.
This research is supported in part
by NSF grant DMS-0070800.

\renewcommand{\theequation}{\thesection.\arabic{equation}}
\renewcommand{\thethm}{\thesection.\arabic{thm}}
\setcounter{equation}{0}
\setcounter{thm}{0}

\section{Differential equations}

Let $V=\coprod_{n\in\mathbb{Z}}V_{(n)}$ be a vertex operator algebra
in the sense of \cite{FLM} and \cite{FHL}.
For a $V$-module $W$, the vertex operator map is a linear map
given by
\begin{eqnarray*}
Y: V&\to& (\mbox{\rm End}\; W)[[x, x^{-1}]]\nn
u&\mapsto &Y(u, x)=\sum_{n\in \Z}u_{n}x^{-n-1},
\end{eqnarray*}
where $x$ is a formal variable and
$u_{n}\in \mbox{\rm End}\; W$ for $n\in \Z$. In this paper, we shall
use
$x, x_{1}, x_{2}, ...$ to denote commuting formal variables and
$z, z_{1}, z_{2}, ...$ to denote complex variables or complex numbers.

For any $V$-module $W$, let $C_{1}(W)$
be the subspace of $W$ spanned by elements of the form $u_{-1}w$ for
$$u\in V_{+}=\coprod_{n>0}V_{(n)}$$
and $w\in W$. If $\dim
W/C_{1}(W)<\infty$, we say that $W$ is {\it $C_{1}$-cofinite} or $W$
satisfies the {\it $C_{1}$-cofiniteness condition}.

In this and the next sections, we shall consider only
$\mathbb{R}$-graded $V$-modules whose subspace spanned by elements of
weights less than or equal to any fixed $r\in \mathbb{R}$ is
finite-dimensional. We shall call such a $V$-module a {\it discretely
$\R$-graded $V$-module}. Clearly, finite direct sums of
irreducible $\R$-graded $V$-modules are discretely $\R$-graded
$V$-modules.

Let $R=\mathbb{C}[z_{1}^{\pm 1}, z_{2}^{\pm 1}, (z_{1}-z_{2})^{-1}]$,
$W_{i}$ ($i=0, 1, 2, 3$) discretely $\R$-graded
$V$-modules satisfying the $C_{1}$-cofiniteness
condition and
$T=R\otimes W_{0}\otimes W_{1}\otimes W_{2}\otimes W_{3}$
which has a natural $R$-module structure. For simplicity, we
shall omit one tensor symbol to write $f(z_{1}, z_{2})\otimes
w_{0}\otimes w_{1}\otimes w_{2}\otimes w_{3}$ as
$f(z_{1}, z_{2})
w_{0}\otimes w_{1}\otimes w_{2}\otimes w_{3}$ in $T$.
For a discretely $\R$-graded $V$-module $W$, as in \cite{FHL}, we use
$W'$ to denote the contragredient
module of $W$. In particular, for $u\in V$
and $n\in \mathbb{Z}$, we have the operators $u_{n}$
on $W'$. For $u\in V$ and $n\in \mathbb{Z}$,
let $u^{*}_{n}:W\to W$ be the adjoint of $u_{n}: W'\to
W'$. Note that since $\wt u_{n}=\wt u-n-1$, we have
$\wt u^{*}_{n}=-\wt u+n+1$.

For $u\in V_{+}$ and $w_{i}\in W_{i}$, $i=0, 1, 2, 3$,
let $J$ be the submodule of $T$ generated by
elements of the form (note the different meanings of the subscripts)
\begin{eqnarray*}
\lefteqn{\mathcal{A}(u, w_{0}, w_{1}, w_{2}, w_{3})}\nn
&&=\sum_{k\ge 0}{-1\choose k}(-z_{1})^{k} u^{*}_{-1-k}
w_{0}\otimes w_{1}\otimes w_{2}\otimes w_{3}
-w_{0}\otimes u_{-1}w_{1}\otimes w_{2}\otimes w_{3}\nn
&& \quad -\sum_{k\ge 0}{-1\choose k}
(-(z_{1}-z_{2}))^{-1-k} w_{0}\otimes
w_{1}\otimes u_{k}w_{2}\otimes w_{3}\nn
&&\quad -\sum_{k\ge 0}{-1\choose k}(-z_{1})^{-1-k}
w_{0}\otimes w_{1}\otimes w_{2}\otimes u_{k}w_{3},
\end{eqnarray*}
\begin{eqnarray*}
\lefteqn{\mathcal{B}(u, w_{0}, w_{1}, w_{2}, w_{3})}\nn
&&=\sum_{k\ge 0}{-1\choose k}(-z_{2})^{k} u^{*}_{-1-k}
w_{0}\otimes w_{1}\otimes w_{2}\otimes w_{3}\nn
&&\quad -\sum_{k\ge 0}{-1\choose k}
(z_{1}-z_{2})^{-1-k} w_{0}\otimes
u_{k}w_{1}\otimes w_{2}\otimes w_{3}
-w_{0}\otimes w_{1}\otimes u_{-1}w_{2}\otimes w_{3}\nn
&&\quad -\sum_{k\ge 0}{-1\choose k}(-z_{2})^{-1-k}
w_{0}\otimes w_{1}\otimes w_{2}\otimes u_{k}w_{3},
\end{eqnarray*}
\begin{eqnarray*}
\lefteqn{\mathcal{C}(u, w_{0}, w_{1}, w_{2}, w_{3})}\nn
&&=u^{*}_{-1}w_{0}\otimes w_{1}\otimes w_{2}\otimes w_{3}
-\sum_{k\ge 0}{-1\choose k}z_{1}^{-1-k}w_{0}\otimes
u_{k}w_{1}\otimes w_{2}\otimes w_{3}\nn
&&\quad -\sum_{k\ge 0}{-1\choose k}z_{2}^{-1-k}w_{0}\otimes
w_{1}\otimes u_{k}w_{2}\otimes w_{3}
-w_{0}\otimes w_{1}\otimes w_{2}\otimes u_{-1}w_{3}
\end{eqnarray*}
and
\begin{eqnarray*}
\lefteqn{\mathcal{D}(u, w_{0}, w_{1}, w_{2}, w_{3})}\nn
&&=u_{-1}w_{0}\otimes w_{1}\otimes w_{2}\otimes w_{3}\nn
&&\quad -\sum_{k\ge 0}{-1\choose k}z_{1}^{1+k}w_{0}\otimes
e^{z_{1}^{-1}L(1)}(-z_{1}^{2})^{L(0)}u_{k}(-z_{1}^{-2})^{L(0)}
e^{-z_{1}^{-1}L(1)}
w_{1}\otimes w_{2}\otimes w_{3}\nn
&& \quad -\sum_{k\ge 0}{-1\choose k}z_{2}^{1+k}w_{0}\otimes
w_{1}\otimes e^{z^{-1}_{2}L(1)}(-z_{2}^{2})^{L(0)}
u_{k}(-z_{2}^{-2})^{L(0)}e^{-z^{-1}_{2}L(1)}w_{2}\otimes w_{3}\nn
&&\quad -w_{0}\otimes w_{1}\otimes w_{2}\otimes u^{*}_{-1}w_{3}.
\end{eqnarray*}

The gradings on $W_{i}$ for $i=0, 1, 2, 3$ induce a
grading (also called {\it weight}) on $W_{0}\otimes
W_{1}\otimes W_{2}\otimes W_{3}$ and then also
on $T$ (here we define the weight of elements of $R$ to be $0$).
Let $T_{(r)}$  be the
homogeneous subspace of weight $r$ for $r\in \mathbb{R}$.
Then $T=\coprod_{r\in \mathbb{R}}T_{(r)}$, $T_{(r)}$
for $r\in \mathbb{R}$ are finitely-generated $R$-modules
and $T_{(r)}=0$ when  $r$ is
sufficiently small. We also let $F_{r}(T)=\coprod_{s\le r}T_{(s)}$
for $r\in \mathbb{R}$.
Then $F_{r}(T)$, $r\in \mathbb{R}$, are finitely-generated
$R$-modules,
$F_{r}(T)\subset F_{s}(T)$ for $r\le s$ and
$\cup_{r\in \mathbb{R}}F_{r}(T)=T$. Let $F_{r}(J)=J\cap F_{r}(T)$
for $r\in \mathbb{R}$. Then $F_{r}(J)$ for $r\in \mathbb{R}$ are
finitely-generated $R$-modules,
$F_{r}(J)\subset F_{s}(J)$ for $r\le s$ and
$\cup_{r\in \mathbb{R}}F_{r}(J)=J$.

\begin{prop}\label{decomposition}
There exists $M\in \mathbb{Z}$ such that for any $r\in \mathbb{R}$,
$F_{r}(T)\subset F_{r}(J)+F_{M}(T)$. In particular, $T=J+F_{M}(T)$.
\end{prop}
\pf
Since $\dim W_{i}/C_{1}(W_{i})<\infty$ for $i=0, 1, 2, 3$,
there exists $M\in \mathbb{Z}$ such that
\begin{eqnarray}\label{subset}
\lefteqn{\coprod_{n>M}T_{(n)}\subset R (C_{1}(W_{0})\otimes W_{1}
\otimes W_{2}\otimes W_{3}) + R (W_{0}\otimes C_{1}(W_{1})\otimes
W_{2}\otimes W_{3})}\nn
&&\quad+ R (W_{0}\otimes W_{1}\otimes
C_{1}(W_{2})\otimes W_{3})+R (W_{0}\otimes W_{1}\otimes
W_{2}\otimes C_{1}(W_{3})).\nn
\end{eqnarray}

We use induction on $r\in \mathbb{R}$.
If $r$ is equal to $M$,
$F_{M}(T)\subset F_{M}(J)+F_{M}(T)$. Now we assume that
$F_{r}(T)\subset F_{r}(J)+F_{M}(T)$ for $r< s$ where $s>M$.
We want to show that any homogeneous element of
$T_{(s)}$ can  be written as
a sum of an element of $F_{s}(J)$ and an element of $F_{M}(T)$.
Since $s>M$, by (\ref{subset}), any element of $T_{(s)}$ is an element
of
the right-hand side of (\ref{subset}). We shall discuss only
the case that this element is
in $R (W_{0}\otimes C_{1}(W_{1})\otimes
W_{2}\otimes W_{3})$; the other cases are completely similar.

We need only discuss elements of the form
$w_{0}\otimes u_{-1}w_{1}\otimes w_{2}\otimes w_{3}$ where
$w_{i}\in W_{i}$ for $i=0, 1, 2, 3$ and $u\in V_{+}$.
By assumption, the weight of
$w_{0}\otimes u_{-1}w_{1}\otimes w_{2}\otimes w_{3}$ is $s$
and the weights of $u^{*}_{-1-k}
w_{0}\otimes w_{1}\otimes w_{2}\otimes w_{3}$,
$w_{0}\otimes w_{1}\otimes w_{2}\otimes u_{k}w_{3}$ and
$w_{0}\otimes
w_{1}\otimes u_{k}w_{2}\otimes w_{3}$ for $k\ge 0$,
are all less than the weight of
$w_{0}\otimes u_{-1}w_{1}\otimes w_{2}\otimes w_{3}$.
So $\mathcal{A}(u, w_{0}, w_{1}, w_{2}, w_{3})
\in F_{s}(J)$.
Thus we see that $w_{0}\otimes u_{-1}w_{1}\otimes w_{2}\otimes w_{3}$
can be written as a sum of an element of $F_{s}(J)$ and
elements of $T$ of weights less than $s$. By the induction assumption,
we know that $w_{0}\otimes u_{-1}w_{1}\otimes w_{2}\otimes w_{3}$
can be written as a sum of an element of $F_{s}(J)$ and
an element of $F_{M}(T)$.

Now we have
\begin{eqnarray*}
T&=&\cup_{r\in \mathbb{R}}F_{r}(T)\nn
&\subset &\cup_{r\in \mathbb{R}}
F_{r}(J)+F_{M}(T)\nn
&=&J+F_{M}(T).
\end{eqnarray*}
But we know that
$J+F_{M}(T)\subset T$. Thus we have $T=J+F_{M}(T)$.
\epfv

We immediately obtain the following:

\begin{cor}
The quotient $R$-module $T/J$ is finitely generated.
\end{cor}
\pf
Since $T=J+F_{M}(T)$ and $F_{M}(T)$
is finitely-generated, $T/J$ is finitely-generated.
\epfv

For an element $\mathcal{W}\in T$, we shall use $[\mathcal{W}]$
to denote the equivalence class in $T/J$ containing $\mathcal{W}$.
We also have:

\begin{cor}
For any $w_{i}\in W_{i}$ ($ i=0, 1, 2, 3$), let $M_{1}$ and $M_{2}$
be the $R$-submodules of $T/J$ generated by $[w_{0}\otimes
L(-1)^{j}w_{1}
\otimes w_{2}\otimes w_{3}]$, $j\ge 0$, and by $[w_{0}\otimes w_{1}
\otimes L(-1)^{j}w_{2}\otimes w_{3}]$, $j\ge 0$, respectively.
Then $M_{1}$ and $M_{2}$ are finitely generated. In particular,
for any $w_{i}\in W_{i}$ ($i=0, 1, 2, 3$), there exist
$a_{k}(z_{1}, z_{2}),
b_{l}(z_{1}, z_{2})\in R$
for $k=1, \dots, m$ and $l=1, \dots, n$ such that
\begin{eqnarray}
&[w_{0}\otimes L(-1)^{m}w_{1}\otimes w_{2}\otimes w_{3}]
+ a_{1}(z_{1}, z_{2})
[w_{0}\otimes L(-1)^{m-1}w_{1}\otimes w_{2}\otimes w_{3}]&\nn
&+\cdots +a_{m}(z_{1}, z_{2})
[w_{0}\otimes w_{1}\otimes w_{2}\otimes
w_{3}]=0,&\label{dependence1}\\
&[w_{0}\otimes w_{1}\otimes L(-1)^{n}w_{2}\otimes w_{3}]
+ b_{1}(z_{1}, z_{2})
[w_{0}\otimes w_{1}\otimes L(-1)^{n-1}w_{2}\otimes w_{3}]&\nn
&+\cdots +b_{n}(z_{1}, z_{2})
[w_{0}\otimes w_{1}\otimes w_{2}\otimes w_{3}]=0.&\label{dependence2}
\end{eqnarray}
\end{cor}
\pf
Since $R$ is a Noetherian ring, any $R$-submodule of the
finitely-generated
$R$-module $T/J$ is also finitely generated. In particular,
$M_{1}$ and $M_{2}$ are finitely generated. The second conclusion
follows immediately.
\epfv

Now we establish the existence of systems of differential equations:

\begin{thm}\label{sys1}
Let $W_{i}$ for $i=0, 1, 2, 3$ be discretely $\R$-graded
$V$-modules satisfying the
$C_{1}$-cofiniteness condition. Then
for any $w_{i}\in W_{i}$ ($i=0, 1, 2, 3$), there exist
$$a_{k}(z_{1}, z_{2}),
b_{l}(z_{1}, z_{2})\in
\mathbb{C}[z_{1}^{\pm 1}, z_{2}^{\pm 1}, (z_{1}-z_{2})^{-1}]$$
for $k=1, \dots, m$ and $l=1, \dots, n$ such that
for any discretely $\R$-graded $V$-modules $W_{4}$, $W_{5}$ and
$W_{6}$, any
intertwining operators $\mathcal{Y}_{1}$, $\mathcal{Y}_{2}$,
$\mathcal{Y}_{3}$, $\mathcal{Y}_{4}$, $\mathcal{Y}_{5}$
and $\mathcal{Y}_{6}$ of types
${W'_{0}\choose W_{1}W_{4}}$, ${W_{4}\choose W_{2}W_{3}}$,
${W_{5}\choose W_{1}W_{2}}$, ${W'_{0}\choose W_{5}W_{3}}$,
${W'_{0}\choose W_{2}W_{6}}$ and ${W_{6}\choose W_{1}W_{3}}$,
respectively, the series
\begin{equation}\label{prod1}
\langle w_{0}, \mathcal{Y}_{1}(w_{1}, z_{1})
\mathcal{Y}_{2}(w_{2}, z_{2})w_{3}\rangle,
\end{equation}
\begin{equation}\label{iter}
\langle w_{0}, \mathcal{Y}_{4}(\mathcal{Y}_{3}(w_{1}, z_{1}-z_{2})
w_{2}, z_{2})w_{3}\rangle
\end{equation}
and
\begin{equation}\label{prod2}
\langle w_{0}, \mathcal{Y}_{5}(w_{2}, z_{2})
\mathcal{Y}_{6}(w_{1}, z_{1})w_{3}\rangle,
\end{equation}
satisfy the expansions of the system of differential equations
\begin{eqnarray}
\frac{\partial^{m}\varphi}{\partial z_{1}^{m}}+a_{1}(z_{1}, z_{2})
\frac{\partial^{m-1}\varphi}{\partial z_{1}^{m-1}}
+\cdots + a_{m}(z_{1}, z_{2})\varphi&=&0,\label{eqn1}\\
\frac{\partial^{n}\varphi}{\partial z_{2}^{n}}+b_{1}(z_{1}, z_{2})
\frac{\partial^{n-1}\varphi}{\partial z_{2}^{n-1}}
+\cdots + b_{n}(z_{1}, z_{2})\varphi&=&0\label{eqn2}
\end{eqnarray}
in the regions $|z_{1}|>|z_{2}|>0$, $|z_{2}|>|z_{1}-z_{2}|>0$
and $|z_{2}|>|z_{1}|>0$,
respectively.
\end{thm}
\pf
For simplicity, we assume that $w_{i}\in W_{i}$ for $i=0, 1, 2, 3$
are homogeneous. Let $\Delta=\wt w_{0}-\wt w_{1}-\wt w_{2}-\wt w_{3}$.
Let $\C(\{x\})$ be the space of all
series of the form $\sum_{n\in \R}a_{n}x^{n}$, $a_{n}\in \C$
for $n\in \R$, such that $a_{n}=0$ when the real part
of $n$ is sufficiently negative. For any complex variable $z$,
we can substitute $z^{n}=e^{n\log z}$
for $x^{n}$ for $n\in \R$ to
obtain a space $\C(\{z\})$ isomorphic to $\C(\{x\})$ (here and below
we use
the convention that $\log z$ is the value of the logarithm of $z$
such that $0\le \arg z<2\pi$). For  formal variables $x_{1}$ and
$x_{2}$,
we consider the
space $\C(\{z_{2}/z_{1}\})[x_{1}^{\pm 1}, x_{2}^{\pm 1}]$. Let
$I_{z_{1}, z_{2}}$ be the subspace of
$\C(\{z_{2}/z_{1}\})[x_{1}^{\pm 1}, x_{2}^{\pm 1}]$ spanned by
elements of the form
$$\left(\sum_{n\in \R}a_{n}(z_{2}/z_{1})^{n}\right)x_{1}^{k}x_{2}^{l}-
\left(\sum_{n\in
\R}a_{n}(z_{2}/z_{1})^{n+i}\right)x_{1}^{k+i}x_{2}^{l-i}$$
for $k, l, i\in \Z$. Let
$$\C(\{z_{2}/z_{1}\})[z_{1}^{\pm 1}, z_{2}^{\pm 1}]
=\C(\{z_{2}/z_{1}\})[x_{1}^{\pm 1}, x_{2}^{\pm 1}]/I_{z_{1}, z_{2}}.$$
Then clearly
$\C(\{z_{2}/z_{1}\})[z_{1}^{\pm 1}, z_{2}^{\pm 1}]$ is a
$\mathbb{C}[[z_{2}/z_{1}]][z_{1}^{\pm 1}, z_{2}^{\pm 1}]$-module and
can be identified with subspaces of $\C\{z_{1}, z_{2}\}$ which
are the space of all series of the form $\sum_{m, n\in \R}b_{m, n}
e^{m\log z_{1}}e^{n\log z_{2}}$, $b_{m, n}\in \C$ for $m, n\in \R$.
Similarly we
have $\mathbb{C}(\{(z_{1}-z_{2})/z_{1}\})[z_{2}^{\pm 1},
(z_{1}-z_{2})^{\pm 1}]$ and
$\C(\{z_{1}/z_{2}\})[z_{1}^{\pm 1}, z_{2}^{\pm 1}]$, where in
$\mathbb{C}(\{(z_{1}-z_{2})/z_{1}\})[z_{2}^{\pm 1},
(z_{1}-z_{2})^{\pm 1}]$, $z_{2}$ and
$z_{1}-z_{2}$ are viewed as independent complex variables.
The spaces $\mathbb{C}(\{(z_{1}-z_{2})/z_{1}\})[z_{2}^{\pm 1},
(z_{1}-z_{2})^{\pm 1}]$ and
$\C(\{z_{1}/z_{2}\})[z_{1}^{\pm 1}, z_{2}^{\pm 1}]$
are $\mathbb{C}[[(z_{1}-z_{2})/z_{2}]]
[z_{2}^{\pm 1},
(z_{1}-z_{2})^{\pm 1}]$- and
$\mathbb{C}[[z_{1}/z_{2}]][z_{1}^{\pm 1}, z_{2}^{\pm 1}]$-modules,
respectively,
and can also be identified with a subspace of $\C\{z_{1}, z_{2}\}$.
For any $r\in \R$, $z_{1}^{r}\C(\{z_{2}/z_{1}\})[z_{1}^{\pm 1},
z_{2}^{\pm 1}]$,
$z_{2}^{r}\mathbb{C}(\{(z_{1}-z_{2})/z_{1}\})[z_{2}^{\pm 1},
(z_{1}-z_{2})^{\pm 1}]$ and
$z_{2}^{r}\C(\{z_{1}/z_{2}\})[z_{1}^{\pm 1}, z_{2}^{\pm 1}]$ are
also $\mathbb{C}[[z_{2}/z_{1}]][z_{1}^{\pm 1}, z_{2}^{\pm 1}]$-,
$\mathbb{C}[[(z_{1}-z_{2})/z_{2}]]
[z_{2}^{\pm 1},
(z_{1}-z_{2})^{\pm 1}]$- and
$\mathbb{C}[[z_{1}/z_{2}]][z_{1}^{\pm 1}, z_{2}^{\pm 1}]$-modules,
respectively, and are also subspaces of $\C\{z_{1}, z_{2}\}$.

For any discretely $\R$-graded $V$-modules $W_{4}$,  $W_{5}$ and
$W_{6}$, any
intertwining operators $\mathcal{Y}_{1}$, $\mathcal{Y}_{2}$,
$\mathcal{Y}_{3}$, $\mathcal{Y}_{4}$, $\mathcal{Y}_{5}$
and $\mathcal{Y}_{6}$ of types
${W'_{0}\choose W_{1}W_{4}}$, ${W_{4}\choose W_{2}W_{3}}$,
${W_{5}\choose W_{1}W_{2}}$, ${W'_{0}\choose W_{5}W_{3}}$,
${W'_{0}\choose W_{2}W_{6}}$ and ${W_{6}\choose W_{1}W_{3}}$,
respectively, consider the maps
\begin{eqnarray*}
\phi_{\mathcal{Y}_{1}, \mathcal{Y}_{2}}: T&\to& z_{1}^{\Delta}
\C(\{z_{2}/z_{1}\})[z_{1}^{\pm 1}, z_{2}^{\pm 1}],\\
\psi_{\mathcal{Y}_{3}, \mathcal{Y}_{4}}: T&\to& z_{2}^{\Delta}
\mathbb{C}(\{(z_{1}-z_{2})/z_{1}\})[z_{2}^{\pm 1},
(z_{1}-z_{2})^{\pm 1}],\\
\xi_{\mathcal{Y}_{5}, \mathcal{Y}_{6}}: T&\to& z_{2}^{\Delta}
\C(\{z_{1}/z_{2}\})[z_{1}^{\pm 1}, z_{2}^{\pm 1}],
\end{eqnarray*}
defined by
\begin{eqnarray*}
\lefteqn{\phi_{\mathcal{Y}_{1}, \mathcal{Y}_{2}}(f(z_{1}, z_{2})
w_{0}\otimes w_{1}\otimes
w_{2}\otimes w_{3})}\nn
&&=\iota_{|z_{1}|>|z_{2}|>0}(f(z_{1}, z_{2}))\langle w_{0},
\mathcal{Y}_{1}(w_{1},
z_{1})\mathcal{Y}_{2}(w_{2}, z_{2})w_{3}\rangle,
\end{eqnarray*}
\begin{eqnarray*}
\lefteqn{\psi(f(z_{1}, z_{2}) w_{0}\otimes w_{1}\otimes
w_{2}\otimes w_{3})}\nn
&&=\iota_{|z_{2}|>|z_{1}-z_{2}|>0}(f(z_{1}, z_{2}))\langle w_{0},
\mathcal{Y}_{4}(\mathcal{Y}_{3}(w_{1},
z_{1}-z_{2})w_{2}, z_{2})w_{3}\rangle,
\end{eqnarray*}
\begin{eqnarray*}
\lefteqn{\xi_{\mathcal{Y}_{5}, \mathcal{Y}_{6}}(f(z_{1}, z_{2})
w_{0}\otimes w_{1}\otimes
w_{2}\otimes w_{3})}\nn
&&=\iota_{|z_{2}|>|z_{1}|>0}(f(z_{1}, z_{2}))\langle w_{0},
\mathcal{Y}_{5}(w_{2},
z_{2})\mathcal{Y}_{6}(w_{1}, z_{1})w_{3}\rangle,
\end{eqnarray*}
respectively, where
\begin{eqnarray*}
\iota_{|z_{1}|>|z_{2}|>0}: R&\to &
\mathbb{C}[[z_{2}/z_{1}]][z_{1}^{\pm 1}, z_{2}^{\pm 1}],\\
\iota_{|z_{2}|>|z_{1}-z_{2}|>0}: R&\to&
\mathbb{C}[[(z_{1}-z_{2})/z_{2}]]
[z_{2}^{\pm 1},
(z_{1}-z_{2})^{\pm 1}],\\
\iota_{|z_{2}|>|z_{1}|>0}: R&\to &
\mathbb{C}[[z_{1}/z_{2}]][z_{1}^{\pm 1}, z_{2}^{\pm 1}],
\end{eqnarray*}
are the maps  expanding elements
of $R$ as series in the regions $|z_{1}|>|z_{2}|>0$,
$|z_{2}|>|z_{1}-z_{2}|>0$, $|z_{2}|>|z_{1}|>0$, respectively.

Using the definition of $J$ and  the Jacobi identity defining
intertwining operators,
we have $\phi_{\mathcal{Y}_{1}, \mathcal{Y}_{2}}(J)
=\psi_{\mathcal{Y}_{3}, \mathcal{Y}_{4}}(J)
=\xi_{\mathcal{Y}_{5}, \mathcal{Y}_{6}}(J)=0$. (In fact, we
purposely choose
$\mathcal{A}(u, w_{0}, w_{1}, w_{2}, w_{3})$,
$\mathcal{B}(u, w_{0}, w_{1}, w_{2}, w_{3})$,
$\mathcal{C}(u, w_{0}, w_{1}, w_{2}, w_{3})$ and
$\mathcal{D}(u, w_{0}, w_{1}, w_{2}, w_{3})$
to be elements of the intersection of
the right-hand side of (\ref{subset})
and  the kernels of $\phi_{\mathcal{Y}_{1},
\mathcal{Y}_{2}}$, $\psi_{\mathcal{Y}_{3}, \mathcal{Y}_{4}}$ and
$\xi_{\mathcal{Y}_{5}, \mathcal{Y}_{6}}$ for all
$\mathcal{Y}_{i}$, $i=1, \dots, 6$.)
Thus
we have the induced maps
\begin{eqnarray*}
\bar{\phi}_{\mathcal{Y}_{1}, \mathcal{Y}_{2}}: T/J&\to& z_{1}^{\Delta}
\mathbb{C}(\{z_{2}/z_{1}\})[z_{1}^{\pm 1}, z_{2}^{\pm 1}],\\
\bar{\psi}_{\mathcal{Y}_{3}, \mathcal{Y}_{4}}: T/J&\to& z_{2}^{\Delta}
\mathbb{C}(\{(z_{1}-z_{2})/z_{1}\})[z_{2}^{\pm 1},
(z_{1}-z_{2})^{\pm 1}],\\
\bar{\xi}_{\mathcal{Y}_{5}, \mathcal{Y}_{6}}: T/J&\to& z_{2}^{\Delta}
\mathbb{C}(\{z_{1}/z_{2}\})[z_{1}^{\pm 1}, z_{2}^{\pm 1}].
\end{eqnarray*}
Applying $\bar{\phi}_{\mathcal{Y}_{1}, \mathcal{Y}_{2}}$,
$\bar{\psi}_{\mathcal{Y}_{3}, \mathcal{Y}_{4}}$ and
$\bar{\xi}_{\mathcal{Y}_{5}, \mathcal{Y}_{6}}$
to (\ref{dependence1}) and
(\ref{dependence2})
and then use the $L(-1)$-derivative and $L(-1)$-bracket
properties for intertwining operators,
we see that (\ref{prod1}), (\ref{iter}) and (\ref{prod2})
indeed satisfy
the expansions of the system in the regions $|z_{1}>|z_{2}|>0$,
$|z_{2}|>|z_{1}-z_{2}|>0$ and $|z_{2}>|z_{1}|>0$, respectively.
\epfv

\begin{rema}
{\rm Note that in the theorems above,
$a_{k}(z_{1}, z_{2})$ for $k=1, \dots, m-1$ and
$b_{l}(z_{1}, z_{2})$ for $l=1, \dots, n-1$,
and consequently the corresponding
system, are independent of $\mathcal{Y}_{1}$,
$\mathcal{Y}_{2}$, $\mathcal{Y}_{3}$, $\mathcal{Y}_{4}$,
$\mathcal{Y}_{5}$ and $\mathcal{Y}_{6}$.}
\end{rema}

The following result can be proved using the same method and so the
proof is omitted:

\begin{thm}\label{sys2}
Let $W_{i}$ for $i=0, \dots, p+1$ be discretely $\R$-graded
$V$-modules satisfying the
$C_{1}$-cofiniteness condition. Then
for any $w_{i}\in W_{i}$ for $i=0, \dots, p+1$, there exist
\begin{equation}\label{coeff}
a_{k_{l}, \;l}(z_{1}, \dots,  z_{p})\in
\mathbb{C}[z_{1}^{\pm 1}, \dots, z_{p}^{\pm 1},
(z_{1}-z_{2})^{-1}, (z_{1}-z_{3})^{-1}, \dots, (z_{p-1}-z_{p})^{-1}],
\end{equation}
for $k_{l}=1, \dots, m_{l}$ and
$l=1, \dots, p,$
such that for any discretely $\R$-graded $V$-modules $\tilde{W}_{q}$
for $q=1, \dots, p-1$, any
intertwining operators $\mathcal{Y}_{1}, \mathcal{Y}_{2},
\dots, \mathcal{Y}_{p-1}, \mathcal{Y}_{p}$,
of types
${W'_{0}\choose W_{1}\tilde{W}_{1}}$, ${\tilde{W}_{1}\choose
W_{2}\tilde{W}_{2}}, \dots,
{\tilde{W}_{p-2}\choose W_{p-1}\tilde{W}_{p-1}}$,
${\tilde{W}_{p-1}\choose W_{p}W_{p+1}}$,
respectively, the series
\begin{equation}\label{p-prod}
\langle w_{0}, \mathcal{Y}_{1}(w_{1}, z_{1})\cdots
\mathcal{Y}_{p}(w_{p}, z_{p})w_{p+1}\rangle
\end{equation}
satisfy the expansions of the system
of differential equations
$$\frac{\partial^{m_{l}}\varphi}{\partial z_{l}^{m_{l}}}+
\sum_{k_{l}=1}^{m_{l}}a_{k_{l}, \;l}(z_{1}, \dots,  z_{p})
\frac{\partial^{m_{l}-k_{l}}\varphi}{\partial
z_{l}^{m_{l}-k_{l}}}=0,\;\;\;
l=1, \dots, p$$
in the region $|z_{1}|>\cdots |z_{p}|>0$.
\end{thm}

\begin{rema}\label{reduction}
{\rm In this and the next section, we
assume for simplicity that all the discretely $\R$-graded $V$-modules
involved satisfy the
$C_{1}$-cofiniteness condition.  But the conclusions in
Theorems \ref{sys1} and \ref{sys2}
and in theorems in the next section still hold if for one of the
discretely $\R$-graded $V$-modules, the $C_{1}$-cofiniteness condition
is replaced by
a much weaker condition that it is generated by
lowest weight vectors (a lowest weight vector is a vector
$w$ such that $u_{n}w=0$ when
$\wt u-n-1<0$). In fact, using the Jacobi
identity, we can always reduce the proof of the general case to the
case
that the vector in this non-$C_{1}$-cofinite discretely $\R$-graded
$V$-module is a
lowest-weight vector. Then from the explicit expressions of
$\mathcal{A}(u,
w_{0}, w_{1}, w_{2}, w_{3})$, $\mathcal{B}(u, w_{0}, w_{1}, w_{2},
w_{3})$, $\mathcal{C}(u, w_{0}, w_{1}, w_{2}, w_{3})$ and
$\mathcal{D}(u, w_{0}, w_{1}, w_{2}, w_{3})$,  it is clear that every
argument in this and the next section still works if we replace this
non-$C_{1}$-cofinite
discretely $\R$-graded $V$-module by the space of all its
lowest-weight vectors.}
\end{rema}

\renewcommand{\theequation}{\thesection.\arabic{equation}}
\renewcommand{\thethm}{\thesection.\arabic{thm}}
\setcounter{equation}{0}
\setcounter{thm}{0}

\section{The regularity of the singular points}

We need  certain  filtrations on $R$
and on the $R$-module $T$.

We discuss only the filtrations associated to the singular point
$z_{1}=z_{2}$.  For $n\in \mathbb{Z}_{+}$, let
$F_{n}^{(z_{1}=z_{2})}(R)$
 be the vector subspace of $R$ spanned by elements of the form
$f(z_{1},
z_{2})(z_{1}-z_{2})^{-n}$ for $f(z_{1}, z_{2})\in
\mathbb{C}[z_{1}^{\pm}, z_{2}^{\pm}]$.  Then with respect to this
filtration, $R$ is a filtered algebra, that is,
$F_{m}^{(z_{1}=z_{2})}(R)\subset F_{n}^{(z_{1}=z_{2})}(R)$ for $m\le
n$,
$R=\cup_{n\in \mathbb{Z}}F_{n}^{(z_{1}=z_{2})}(R)$ and
$F_{m}^{(z_{1}=z_{2})}(R)F_{n}^{(z_{1}=z_{2})}(R)\subset
F_{m+n}^{(z_{1}=z_{2})}(R)$ for any $m, n\in \mathbb{Z}_{+}$.

For convenience, we shall use $\sigma$ to denote
$\wt w_{0}+\wt w_{1}+\wt w_{2}+\wt w_{3}$ for
$w_{i}\in W_{i}$, $i=0, 1, 2, 3$, when the dependence
on $w_{0}$, $w_{1}$, $w_{2}$ and $w_{3}$ are clear.
Let $F_{r}^{(z_{1}=z_{2})}(T)$ for $r\in \mathbb{R}$
be the subspace of $T$
spanned by elements of the form $f(z_{1}, z_{2})(z_{1}-z_{2})^{-n}
w_{0}
\otimes w_{1}\otimes w _{2}\otimes w_{3}$ where $f(z_{1}, z_{2})\in
\mathbb{C}[z_{1}^{\pm}, z_{2}^{\pm}]$, $n\in \mathbb{Z}_{+}$ and
$w_{i}\in
W_{i}$ ($i=0, 1, 2, 3$) satisfying
$n+\sigma \le r$.
These subspaces give a filtration of $T$ in the following sense:
$F_{r}^{(z_{1}=z_{2})}(T)\subset F_{s}^{(z_{1}=z_{2})}(T)$ for
$r\le s$;
$T=\cup_{r\in \mathbb{R}}F_{r}^{(z_{1}=z_{2})}(T)$;
$F_{n}^{(z_{1}=z_{2})}(R)F_{r}^{(z_{1}=z_{2})}(T)\subset
F_{r+n}^{(z_{1}=z_{2})}(T)$.

Let $F_{r}^{(z_{1}=z_{2})}(J)=F_{r}^{(z_{1}=z_{2})}(T)\cap J$
for $r\in \mathbb{R}$.
We need the following refinement of Proposition \ref{decomposition}:

\begin{prop}\label{filtr-z1=z2}
For any $r\in \mathbb{R}$,
$F_{r}^{(z_{1}=z_{2})}(T)
\subset F_{r}^{(z_{1}=z_{2})}(J)+F_{M}(T)$.
\end{prop}
\pf
The proof is a refinement of the proof of Proposition
\ref{decomposition}. The only additional property we need
is that the elements  $\mathcal{A}(u, w_{0}, w_{1}, w_{2}, w_{3})$,
$\mathcal{B}(u, w_{0}, w_{1}, w_{2}, w_{3})$,
$\mathcal{C}(u, w_{0}, w_{1}, w_{2}, w_{3})$ and
$\mathcal{D}(u, w_{0}, w_{1}, w_{2}, w_{3})$ are all in
$F^{(z_{1}=z_{2})}_{\swt u+\sigma}(J)$.
This  is clear.
\epfv

We also consider the ring $\mathbb{C}[z_{1}^{\pm}, z_{2}^{\pm}]$
and the $\mathbb{C}[z_{1}^{\pm}, z_{2}^{\pm}]$-module
$$T^{(z_{1}=z_{2})}=\mathbb{C}[z_{1}^{\pm}, z_{2}^{\pm}]\otimes
W_{0}\otimes W_{1}\otimes W_{2}\otimes W_{3}.$$
Let $T_{(r)}^{(z_{1}=z_{2})}$ for $r\in \mathbb{R}$ be the space
of elements of $T^{(z_{1}=z_{2})}$ of  weight $r$. Then
$T^{(z_{1}=z_{2})}=\coprod_{r\in \mathbb{R}}T_{(r)}^{(z_{1}=z_{2})}$.

Let $w_{i}\in W_{i}$ for $i=0, 1, 2, 3$. Then by
Proposition \ref{filtr-z1=z2},
$$w_{0}\otimes w_{1}\otimes w_{2}\otimes w_{3}
=\mathcal{W}_{1}+\mathcal{W}_{2}$$
where $\mathcal{W}_{1}\in
F_{\sigma}^{(z_{1}=z_{2})}(J)$
and $\mathcal{W}_{2}\in
F_{M}(T)$.

\begin{lemma}\label{lemma1}
For any $s\in [0, 1)$, there exist
$S\in \mathbb{R}$ such that $s+S\in \mathbb{Z}_{+}$ and for any
$w_{i}\in W_{i}$, $i=0, 1, 2, 3$, satisfying $\sigma\in s+\mathbb{Z}$,
$(z_{1}-z_{2})^{\sigma+S}
\mathcal{W}_{2}\in T^{(z_{1}=z_{2})}$.
\end{lemma}
\pf
Let $S$ be a real number such that $s+S\in \mathbb{Z}_{+}$
and such that for any $r\in \mathbb{R}$ satisfying
$r\le -S$, $T_{(r)}=0$. By definition, elements of
$F_{r}^{(z_{1}=z_{2})}(T)$ for any $r\in \mathbb{R}$
are sums of elements of the form
$f(z_{1}, z_{2})(z_{1}-z_{2})^{-n} \tilde{w}_{0}
\otimes \tilde{w}_{1}\otimes \tilde{w} _{2}\otimes
\tilde{w}_{3}$ where $f(z_{1}, z_{2})\in
\mathbb{C}[z_{1}^{\pm}, z_{2}^{\pm}]$, $n\in \mathbb{Z}_{+}$ and
$\tilde{w}_{i}\in
W_{i}$ ($i=0, 1, 2, 3$) satisfying
$n+\wt \tilde{w}_{0}+\wt \tilde{w}_{1}+ \wt \tilde{w}_{2}
+\wt \tilde{w}_{3}\le r$.
Since $\wt \tilde{w}_{0}+\wt \tilde{w}_{1}+ \wt \tilde{w}_{2}
+\wt \tilde{w}_{3}>-S$, we obtain
$r-n>-S$ or $r+S-n>0$. Thus
$(z_{1}-z_{2})^{r+S}F_{r}^{(z_{1}=z_{2})}(T)\in T^{(z_{1}=z_{2})}$
if $r+S\in \mathbb{Z}$.

By definition,
$$\mathcal{W}_{2}=w_{0}\otimes w_{1}\otimes w_{2}\otimes w_{3}
-\mathcal{W}_{1},$$
where
$$\mathcal{W}_{1}\in
F_{\sigma}^{(z_{1}=z_{2})}(J)
\subset F_{\sigma}^{(z_{1}=z_{2})}(T).$$
By the discussion above,
$(z_{1}-z_{2})^{\sigma+S}
\mathcal{W}_{1}\in T^{(z_{1}=z_{2})}$ and by definition,
$$w_{0}\otimes w_{1}\otimes w_{2}\otimes w_{3}\in T^{(z_{1}=z_{2})}.$$
Thus $(z_{1}-z_{2})^{\sigma+S}
\mathcal{W}_{2}\in T^{(z_{1}=z_{2})}$.
\epfv

\begin{thm}\label{regular}
Let $W_{i}$ and $w_{i}\in W_{i}$ for $i=0, 1, 2, 3$
be the same as in Theorem \ref{sys1}. For any possible singular points
of the form $z_{1}=0$ or $z_{2}=0$ or $z_{1}=\infty$ or $z_{2}=\infty$ or
$z_{1}=z_{2}$ or  $z_{1}^{-1}(z_{1}-z_{2})=0$ or
or $z_{2}^{-1}(z_{1}-z_{2})=0$,
there exist
$$a_{k}(z_{1}, z_{2}),
b_{l}(z_{1}, z_{2})\in
\mathbb{C}[z_{1}^{\pm 1}, z_{2}^{\pm 1}, (z_{1}-z_{2})^{-1}]$$
for $k=1, \dots, m$ and $l=1, \dots, n$, such that
this singular point of the system (\ref{eqn1}) and (\ref{eqn2})
satisfied by (\ref{prod1}), (\ref{iter}) and (\ref{prod2})
is regular.
\end{thm}
\pf
We shall only prove the theorem for the singular points $z_{1}=z_{2}$
and $z_{2}^{-1}(z_{1}-z_{2})=0$. 
The latter case is one of the most interesting case 
since it will give us the expansion of the solutions in the 
region $|z_{2}|>|z_{1}-z_{2}|>0$.  The other cases are similar.

By Proposition \ref{filtr-z1=z2},
$$w_{0}\otimes L(-1)^{k}w_{1}\otimes w_{2}\otimes w_{3}
=\mathcal{W}_{1}^{(k)}+\mathcal{W}_{2}^{(k)}$$
for $k\ge 0$,
where $\mathcal{W}^{(k)}_{1}\in
F_{\sigma+k}^{(z_{1}=z_{2})}(J)$
and $\mathcal{W}^{(k)}_{2}\in
F_{M}(T)$.

By  Lemma \ref{lemma1}, there exists $S\in \mathbb{R}$ such that
$\sigma+S\in \mathbb{Z}_{+}$ and
$$(z_{1}-z_{2})^{\sigma+k+S}
\mathcal{W}^{(k)}_{2}\in T^{(z_{1}=z_{2})}$$
and
thus
$$(z_{1}-z_{2})^{\sigma+k+S}
\mathcal{W}^{(k)}_{2}\in \coprod_{r\le M}T_{(r)}^{(z_{1}=z_{2})}$$
for $k\ge 0$.
Since $\mathbb{C}[z_{1}^{\pm}, z_{2}^{\pm}]$ is a Noetherian ring and
$\coprod_{r\le M}T_{(r)}^{(z_{1}=z_{2})}$ is a finitely-generated
$\mathbb{C}[z_{1}^{\pm}, z_{2}^{\pm}]$-module,
the submodule of $\coprod_{r\le M}T_{(r)}^{(z_{1}=z_{2})}$ generated
by
$(z_{1}-z_{2})^{\sigma+k+S}
\mathcal{W}^{(k)}_{2}$ for $k\ge 0$ is also finitely generated.
Let $(z_{1}-z_{2})^{\sigma+k+S}
\mathcal{W}^{(k)}_{2}$ for $k=0, \dots, m-1$
be a set of generators of this submodule.
Then there exist $c_{k}(z_{1}, z_{2})\in
\mathbb{C}[z_{1}^{\pm}, z_{2}^{\pm}]$
for $k=0,
\dots, m-1$ such that
$$(z_{1}-z_{2})^{\sigma+m+S}
\mathcal{W}^{(m)}_{2}=-\sum_{k=0}^{m-1}c_{k}(z_{1}, z_{2})
(z_{1}-z_{2})^{\sigma+k+S}
\mathcal{W}^{(k)}_{2}$$
or equivalently
$$\mathcal{W}^{(m)}_{2}+\sum_{k=0}^{m-1}c_{k}(z_{1}, z_{2})
(z_{1}-z_{2})^{k-m}
\mathcal{W}^{(k)}_{2}=0.$$
Thus
\begin{eqnarray}\label{w1}
\lefteqn{w_{0}\otimes L(-1)^{m}w_{1}\otimes w_{2}\otimes w_{3}}\nn
&&+\sum_{k=0}^{m-1}c_{k}(z_{1}, z_{2})(z_{1}-z_{2})^{k-m}
w_{0}\otimes L(-1)^{k}w_{1}\otimes w_{2}\otimes w_{3}\nn
&&\quad =\mathcal{W}^{(m)}_{1}+\sum_{k=0}^{m-1}c_{k}(z_{1}, z_{2})
(z_{1}-z_{2})^{k-m}
\mathcal{W}^{(k)}_{1}.
\end{eqnarray}
Since $\mathcal{W}^{(k)}_{1}\in
F_{\sigma+k}^{(z_{1}=z_{2})}(J)
\subset J$, the right-hand side of (\ref{w1})
is in $J$. Thus we obtain
\begin{eqnarray}\label{dependence3}
\lefteqn{[w_{0}\otimes L(-1)^{m}w_{1}\otimes w_{2}\otimes w_{3}]}\nn
&&+\sum_{k=0}^{m-1}c_{k}(z_{1}, z_{2})(z_{1}-z_{2})^{k-m}
[w_{0}\otimes L(-1)^{k}w_{1}\otimes w_{2}\otimes w_{3}]=0.
\end{eqnarray}

Similarly we can find $d_{l}(z_{1}, z_{2})\in
\mathbb{C}[z_{1}^{\pm}, z_{2} ^{\pm}]$
for $l=0, \dots, n$
such that
\begin{eqnarray}\label{dependence4}
\lefteqn{[w_{0}\otimes w_{1}\otimes L(-1)^{n}w_{2}\otimes w_{3}]}\nn
&&+\sum_{l=0}^{n-1}d_{l}(z_{1}, z_{2})(z_{1}-z_{2})^{l-n}
[w_{0}\otimes w_{1}\otimes L(-1)^{l}w_{2}\otimes w_{3}]=0.
\end{eqnarray}
Now it is clear that the singular point $z_{1}=z_{2}$ 
is regular. 

To prove that the singular point $z_{2}^{-1}(z_{1}-z_{2})=0$ 
is also regular, we introduce new gradings on 
$R$ and $T$.
By asigning the degrees of $z_{1}$ and $z_{2}$ to be $-1$, 
we obtain a grading on $R$. Equipped with this grading, $R$
is a $\mathbb{Z}$-graded ring. Together with the grading (by weights)
on $W_{0}\otimes W_{1}\otimes W_{2}\otimes W_{4}$, this grading 
on $R$ gives another grading on $T$ such that $T$ is a graded 
module for $R$. Note that when 
$u, w_{0}, w_{1}, w_{2}, w_{3}$ are homogeneous, the elements
$\mathcal{A}(u, w_{0}, w_{1}, w_{2}, w_{3})$, 
$\mathcal{B}(u, w_{0}, w_{1}, w_{2}, w_{3})$,
$\mathcal{C}(u, w_{0}, w_{1}, w_{2}, w_{3})$ and
$\mathcal{D}(u, w_{0}, w_{1}, w_{2}, w_{3})$ are all 
homogeneous with respect to this new grading. Thus this grading 
induces a grading on $T/J$ such that $T/J$ is also a graded 
$R$-module. Now note that for homogeneous $w_{0}, w_{1}, w_{2}, w_{3}$,
$$[w_{0}\otimes L(-1)^{l}w_{1}\otimes w_{2}\otimes w_{3}],\;\;\;\;
[w_{0}\otimes w_{1}\otimes L(-1)^{l}w_{2}\otimes w_{3}]\in T/J$$
for $l=0, \dots, n$ are  homogeneous of weights 
$\wt w_{0}+\wt w_{1}+\wt w_{2}+\wt w_{3}+l$, respectively. Thus 
by comparing the degrees of the terms in 
(\ref{dependence3}) and (\ref{dependence4}), we see that 
we can always 
choose $c_{k}(z_{1}, z_{2})$ and $d_{l}(z_{1}, z_{2})$ for 
$k=1, \dots, m$, $l=0, \dots, n$ 
to be elements of $R$ of degree $0$. 
Since $c_{k}(z_{1}, z_{2})$ for $k=0, \dots, m$ and 
$d_{l}(z_{1}, z_{2})$ for $l=0, \dots, n$ are of degree
$0$, they are actually functions of 
$z_{3}=z_{2}^{-1}(z_{1}-z_{2})$. Using this fact and 
changing variables from $z_{1}$ and $z_{2}$ to $z_{3}$ and $z_{2}$,
we see immediately
that the singular point $z_{3}=0$ is
regular. 
\epfv

\begin{rema}\label{regular-2-var}
{\rm The regularity of the singular points proved in 
the theorem above is enough for our applications. In fact 
the sets of the singular points discussed in the theorem above are
all given by one equation and thus the proof of the regularity
is  the same as the proof of the regularity of singular points
of ordinary differential equations. Using the 
standard method in the theory of differential 
equations, we can actually prove the regularity of the sets
of singular points given by two equations such as 
$z_{1}=z_{2}$, $z_{2}=\infty$. Though it is not needed 
in the present paper and in our theory, we sketch a proof here.
Let $\hat{\varphi}$ be the vector whose 
components are 
$$\left((z_{1}-z_{2})\frac{\partial}{\partial z_{1}}\right)^{i}
\left((z_{1}-z_{2})\frac{\partial}{\partial z_{2}}\right)^{j}\varphi$$
for $i=0, \dots, m-1$ and $j=0, \dots, n-1$. 
We also change the variables
to $z_{3}=z_{2}^{-1}(z_{1}-z_{2})$ and 
$z_{4}=z^{-1}_{2}$. 
Then the system (\ref{eqn1}) and (\ref{eqn2}) now becomes
\begin{eqnarray}
z_{3}\frac{\partial}{\partial z_{3}}\hat{\varphi}&=&
A(z_{3})\hat{\varphi},\label{c-var1}\\
z_{4}\frac{\partial}{\partial z_{4}}\hat{\varphi}
&=&\frac{-(1+z_{3})A(z_{3})-B(z_{3})}{z_{3}}\hat{\varphi}
\label{c-var2-1}
\end{eqnarray}
Using the theory of ordinary differential equations of 
regular singular points (see, for example, Appendix B
of \cite{K}), we can solve (\ref{c-var1}) first and then 
substitute the solutions into (\ref{c-var2-1}) to 
solve the 
system (\ref{eqn1}) 
and (\ref{eqn2}). From the explicit form of the solutions,
we see immediately that the singular points $z_{3}=0$, $z_{4}=0$ of the 
system (\ref{c-var1}) and
(\ref{c-var2-1}) are regular. So the 
the singular points $z_{3}=0$, $z_{4}=0$ of the 
system (\ref{eqn1}) 
and (\ref{eqn2}) are also regular. 

In the literature, there is actually a stronger definition 
of regularity of singularities in the case of several 
variables. A set of singular points 
are called {\it regular} in this strong sense if the system 
is meromorphically equivalent to a compatible (or integrable) system
with singularities of logarithmic type (see, for example, Chapter 6 of
\cite{BK} for definitions of singularities of logarithmic type and
regular singularities). Let $G(z_{3})$ be 
the $mn\times mn$ invertible matrix 
meromorphic 
in $z_{3}$ such that 
$G(z_{3})\frac{-(1+z_{3})A(z_{3})-B(z_{3})}{z_{3}}G(z_{3})^{-1}$ is a
Jordan canonical form. Use $G(z_{3})$
as the meromorphic gauge transformation. Then
it is easy to see that the system (\ref{c-var1}) and
(\ref{c-var2-1}) becomes a system whose singular points 
$z_{3}=0$, $z_{4}=0$ are of logarithmic type. 
Thus the singular points $z_{3}=0$, $z_{4}=0$ of the
system (\ref{c-var1}) and
(\ref{c-var2-1}) are also regular in this strong sense.}
\end{rema}

The following result is proved in the same way and so the proof
is omitted:

\begin{thm}\label{sys2-regular}
For any possible
singular point
of the system in Theorem \ref{sys2} of the form $z_{i}=0$
or $z_{i}=\infty$ or $z_{i}=z_{j}$ for $i\ne j$, there exist
elements as in (\ref{coeff}) such that this singular point of
the system in Theorem \ref{sys2} is regular.
\end{thm}

\renewcommand{\theequation}{\thesection.\arabic{equation}}
\renewcommand{\thethm}{\thesection.\arabic{thm}}
\setcounter{equation}{0}
\setcounter{thm}{0}

\section{Applications}

In this section, we prove the finiteness of the fusion rules
and the convergence and extension property needed in the 
tensor product theory for the category of modules for a
vertex operator algebra developed by Lepowsky 
and the author in \cite{HL1}--\cite{HL4} and
\cite{H1} and in the theory of intertwtining operator algebras
developed by the author in \cite{H1} and \cite{H3}--\cite{H6}. 
We also discuss conditions which imply some representation-theoretic
conditions needed in these theories. Using all these results,
we obtain our main theorems on the constructions of vertex 
tensor category structures and intertwining operator algebras. 

First we have:

\begin{thm}\label{fusion}
Let $V$ be a vertex operator algebra and $W_{1}$, $W_{2}$
and $W_{3}$ three discretely $\R$-graded $V$-modules. If $W_{1}$,
$W_{2}$
and $W'_{3}$ are $C_{1}$-cofinite, then the
fusion rule among $W_{1}$, $W_{2}$
and $W_{3}$ is finite.
\end{thm}
\pf
Note that $V$ is $C_{1}$-cofinite.
Consider $W_{3}'\otimes V\otimes W_{1}\otimes W_{2}$
and the corresponding $R$-module
$$T=R\otimes W_{3}'\otimes V\otimes W_{1}\otimes W_{2}.$$
We fix
$z^{(0)}_{1}, z^{(0)}_{2}\in \mathbb{C}$
satisfying $z^{(0)}_{1}, z^{(0)}_{2}\ne 0$
and $z^{(0)}_{1}\ne z^{(0)}_{2}$. Consider the evaluation maps {}from
$R$
to $\mathbb{C}$
and {}from $T$ to $W_{3}'\otimes V\otimes W_{1}\otimes W_{2}$ given
by evaluating elements of $R$ and $T$ at $(z_{1}^{(0)}, z_{2}^{(0)})$.
We use $E$ to denote these maps.
Since $T/J$ is a finitely-generated $R$-module, $E(T)/E(J)$ is
a finite-dimensional vector space. For any intertwining
operator $\mathcal{Y}$ of type ${W_{3}\choose W_{1} W_{2}}$,
the map $\bar{\phi}_{Y_{3}, \mathcal{Y}}$ induces a map
$E(\bar{\phi}_{Y_{3}, \mathcal{Y}})$ {}from $E(T)/E(J)$ to
$\mathbb{C}$, where $Y_{3}$ is the vertex operator map
defining the $V$-module structure on $W_{3}$.
Thus we obtain a linear map {}from the
space of intertwining operators of type  ${W_{3}\choose W_{1} W_{2}}$
to the space of linear maps {}from $E(T)/E(J)$ to
$\mathbb{C}$. Since intertwining operators are actually determined
by their values at one point (see \cite{HL1}),
we see that this linear map is injective.
Thus  the dimension of the space of intertwining operators of type
${W_{3}\choose W_{1} W_{2}}$ is less than or equal to the
dimension of the dual space of $E(T)/E(J)$. Since
$E(T)/E(J)$ is finite dimensional, we see that the dimension of its
dual space and therefore the dimension of the space of
intertwining operators of type
${W_{3}\choose W_{1} W_{2}}$ is finite.
\epfv

\begin{rema}
{\rm In fact, the conclusion of Theorem \ref{fusion} still holds
if for one of the
discretely $\R$-graded $V$-modules $W_{1}$, $W_{2}$ and $W_{3}$,
the $C_{1}$-cofiniteness condition is replaced by
the condition that it is generated by
lowest weight vectors. The reason is the same as in
Remark \ref{reduction}. In particular, the fusion rule is finite
when $W_{1}$ and $W_{2}$ are $C_{1}$-cofinite and $W_{3}$ is
completely reducible.}
\end{rema}

\begin{rema}
{\rm Theorem \ref{fusion} has been known to physicists for some time
and
originally goes back to the work of  Nahm \cite{N}. In \cite{L}, Li
proved that for irreducible $V$-modules satisfying a slightly weaker
cofiniteness condition,
the fusion rule is finite. In particular, in the case that
$W_{1}$, $W_{2}$ and $W_{3}$ are irreducible, Theorem \ref{fusion} is
a consequence. In \cite{AN}, Abe and Nagatomo proved that the
dimensions
of the spaces of conformal blocks on the Riemann sphere are finite
if the modules involved are $C_{1}$-cofinite (called $B_{1}$
finite in \cite{AN}) and the vertex operator algebra satisfies some
additional conditions. In particular, under
these conditions, the fusion rules are finite. }
\end{rema}

To formulate the next result,
we need:

\begin{defn}
{\rm Let $V$ be a vertex operator algebra.
We say that  products of intertwining operators for $V$
have the {\it
convergence and extension property}
if for any $\C$-graded $V$-modules $W_{i}$ ($i=0, 1, 2, 3$) and
any intertwining operators
${\cal Y}_{1}$ and ${\cal Y}_{2}$ of types
${W'_{0}}\choose {W_{1}W_{4}}$ and ${W_{4}}\choose
{W_{2}W_{3}}$, respectively, there exists
an integer $N$ (depending only on ${\cal Y}_{1}$ and ${\cal Y}_{2}$),
and for any $w_{i}\in W_{i}$ ($i=0, 1, 2, 3$), there exist
$r_{i}, s_{i}\in \mathbb{R}$ and analytic
functions $f_{i}(z)$ on $|z|<1$ for $i=1, \dots, j$
satisfying
\begin{equation}\label{si}
\wt w_{(1)}+\wt w_{(2)}+s_{i}>N,\;\;\;i=1, \dots, j,
\end{equation}
such that (\ref{prod1}) is absolutely
convergent when $|z_{1}|>|z_{2}|>0$ and can be analytically extended
to
the multivalued analytic function
\begin{equation}\label{phyper}
\sum_{i=1}^{j}z_{2}^{r_{i}}(z_{1}-z_{2})^{s_{i}}
f_{i}\left(\frac{z_{1}-z_{2}}{z_{2}}\right)
\end{equation}
when $|z_{2}|>|z_{1}-z_{2}|>0$.}
\end{defn}

This property was introduced and needed in \cite{H1}
in order to construct the
associativity isomorphism between two iterated
tensor products of three modules for a suitable vertex operator
algebra.
Also recall that a $\mathbb{C}$-graded generalized module for a vertex
operator
algebra is a $\mathbb{C}$-graded vector space equipped
with a vertex operator map satisfying all the
axioms for modules except the two grading-restriction
axioms.

We have:

\begin{thm}\label{conv-ext}
Let $V$ be a vertex operator algebra satisfying the following
conditions:

\begin{enumerate}

\item  Every $\C$-graded generalized $V$-module is a direct sum of
$\C$-graded irreducible
$V$-modules.

\item There are only finitely many inequivalent $\C$-graded
irreducible $V$-modules and they are all $\R$-graded.

\item Every
$\R$-graded irreducible $V$-module satisfies
the $C_{1}$-cofiniteness condition.

\end{enumerate}
Then every $\C$-graded $V$-module is a discretely $\R$-graded
$V$-module
and products of intertwining operators for $V$
have the convergence and extension property. In addition,
for any discretely $\R$-graded $V$-modules (or equivalently,
$\C$-graded $V$-modules), intertwining operators and elements of
the $V$-modules as in Theorem \ref{sys2},
(\ref{p-prod}) is absolutely convergent when $|z_{1}|>\cdots
|z_{p}|>0$ and can be analytically extended to the region 
given by $z_{i}\ne z_{j}$ ($i\ne j$) $z_{i}\ne 0$ ($i=1, \dots, n$).
\end{thm}
\pf
It is clear that every $\C$-graded $V$-module
is a discretely $\R$-graded $V$-module. In the remaining part of
this section, if every $\C$-graded $V$-module
is a discretely $\R$-graded $V$-module, we shall call
discretely $\R$-graded $V$-modules simply $V$-modules.

For any $w_{i}\in W_{i}$, $i=0, 1,2 ,3$, using Theorems \ref{sys1} and
\ref{regular}, the theory of differential equations of regular singular
points (see, for example, Appendix B of \cite{K}), it is easy to show
that (\ref{prod1}) is absolutely convergent when $|z_{1}|>|z_{2}|>0$. In
fact, the expansion coefficients of (\ref{prod1}) as a series in powers
of $z_{2}$ are series in powers of $z_{1}$ of the form $\langle w_{0},
\mathcal{Y}_{1}(w_{1}, z_{1})w_{4}\rangle$ for $w_{4}\in W_{4}$. These
series in powers of $z_{1}$ have only finitely many terms and thus are
convergent for any $z_{1}\ne 0$. So for fixed $z_{1}\ne 0$,
(\ref{prod1}) is a series in powers of $z_{2}$.  Now fix $z_{1}\ne 0$.
By Theorem \ref{regular} the equation (\ref{eqn2}) as an ordinary
differential equation with the variable $z_{2}$ in the region
$|z_{1}|>|z_{2}|>0$ has a regular singular point at $z_{2}=0$. Since
(\ref{prod1}) satisfies (\ref{eqn2}), it must be convergent 
as a series in powers of $z_{2}$. Since the coefficients of 
(\ref{eqn2}) are analytic in $z_{1}$, the sum of (\ref{prod1})
is also analytic in $z_{1}$. So the sum is in fact analytic in both
$z_{1}$ and $z_{2}$ and its expansion as a double series in powers of 
$z_{1}$ and $z_{2}$ must be absolutely convergent. 

To prove the remaining part of the convergence and 
extension property, we need the regularity of the 
system (\ref{eqn1}) 
and (\ref{eqn2}) at the the singular point $z_{2}^{-1}(z_{1}-z_{2})=0$. 
Using the $L(0)$-conjugation property for intertwining operators,
we see that (\ref{prod1}) is equal to 
$$
z_{2}^{\swt w_{0}-\swt w_{2}-\swt w_{2}-\swt w_{2}}\langle w_{0}, 
\mathcal{Y}_{1}(w_{1}, 1+z_{2}^{-1}(z_{1}-z_{2}))
\mathcal{Y}_{2}(w_{2}, 1)w_{3}\rangle.
$$
By Theorems \ref{sys1}, the series 
(\ref{prod1}) satisfies (\ref{eqn1}) and 
by changing the variables to $z_{3}=z_{2}^{-1}(z_{1}-z_{2})$
and $z_{2}$, 
$$\langle w_{0}, 
\mathcal{Y}_{1}(w_{1}, 1+z_{3})
\mathcal{Y}_{2}(w_{2}, 1)w_{3}\rangle$$
satisfies an ordinary differential 
equation with the variable $z_{3}$.
By Theorem \ref{regular}, we see that the singular point 
$z_{3}=0$ of this ordinatry differential
equation is regular. 
By the theory of (ordinary) differential equations 
of regular singular
points (see, for example, Appendix B of \cite{K})
and the fact that there are no
logarithmic terms in the two intertwining operators in (\ref{prod1}),
we see  that
there exist $r_{i, k}, s_{i, k}\in \mathbb{R}$ and analytic
functions $f_{i, k}(z)$ on $|z|<1$ for $i=1, \dots, j$, $k=1,
\dots, K$, such that (\ref{prod1}) can be analytically extended to the
multivalued analytic function
\begin{equation}\label{extension}
\sum_{k=1}^{K}\sum_{i=1}^{j}z_{2}^{r_{i, k}}(z_{1}-z_{2})^{s_{i, k}}
f_{i, k}\left(\frac{z_{1}-z_{2}}{z_{2}}\right)
\left(\log\left(\frac{z_{1}-z_{2}}{z_{2}}\right)\right)^{k}
\end{equation}
in the region $|z_{2}|>|z_{1}-z_{2}|>0$. Since (\ref{prod1})
for general elements $w_{i}\in W_{i}$, $i=0, 1, 2, 3$,
are $R$-linear combinations
of (\ref{prod1}) for those $w_{i}\in W_{i}$ ($i=0, 1, 2, 3$)
satisfying $\wt w_{0}+\wt w_{1}+\wt w_{2}+\wt w_{3}\le M$,
we see that $K$ can be taken to be independent of
$w_{i}\in W_{i}$, $i=0, 1, 2, 3$.
We now show that
$K=0$.

We see that (\ref{extension}) can be written as
\begin{equation}\label{extension2}
\sum_{k=1}^{K}g_{k}(w_{0}, w_{1}, w_{2}, w_{3}; z_{1}, z_{2})
\left(\log\left(1-\frac{z_{2}}{z_{1}}\right)\right)^{k},
\end{equation}
where $g_{k}(w_{0}, w_{1}, w_{2}, w_{3};
z_{1}, z_{2})$ for $k=1, \dots, K$ are linear combinations
of products of $z_{2}^{r_{i, k}}$, $(z_{1}-z_{2})^{s_{i, k}}$,
$f_{i, k}\left(\frac{z_{1}-z_{2}}{z_{2}}\right)$ and
$\left(\log\left(
\frac{z_{1}}{z_{2}}\right)\right)^{k}$ for $i=1, \dots, j$ and
$k=1, \dots, K$. Since (\ref{extension}) and
$\left(\log\left(1-\frac{z_{2}}{z_{1}}\right)\right)^{k}$
for $k=1, \dots, K$ can be analytically extended to
the region $|z_{1}|>|z_{2}|>0$, $g_{k}(w_{0}, w_{1}, w_{2}, w_{3};
z_{1}, z_{2})$ for $k=1, \dots, K$ can also be extended to
this region. Consequently (\ref{prod1}) is equal to
(\ref{extension2}).

Assume that $K\ne 0$. We can view $g_{K}(w_{0}, w_{1}, w_{2}, w_{3};
z_{1}, z_{2})$ as the value at $w_{1}\otimes w_{2}\otimes w_{3}$
of the image of $w_{0}$ under a linear map {}from $W_{0}$ to
$(W_{1}\otimes W_{2}\otimes W_{3})^{*}$.
Then by the
properties of (\ref{prod1}) and the linear independence
of $\left(\log\left(1-\frac{z_{2}}{z_{1}}\right)\right)^{k}$,
$k=1, \dots, K$,  the image of $w_{0}$ under
this linear map {}from $W_{0}$ to
$(W_{1}\otimes W_{2}\otimes W_{3})^{*}$ satisfies the
$P(z_{1}, z_{2})$-compatibility condition,
the $P(z_{1}, z_{2})$-local grading-restriction condition and
the $P(z_{2})$-local grading-restriction condition (see \cite{H1}).
Using Conditions 1--3, Theorem \ref{fusion}
Theorem 14.10 in \cite{H1} and
the construction of the associativity isomorphisms
in \cite{H1}, we see that this image is in fact in
$\Psi_{P(z_{1}, z_{2})}^{(1)}(W_{1}\hboxtr_{P(z_{1})}(W_{1}\boxtimes
W_{3}))$
where
$$\Psi_{P(z_{1}, z_{2})}^{(1)}: W_{1}\hboxtr_{P(z_{1})}(W_{1}\boxtimes
W_{3})
\to (W_{1}\otimes W_{2}\otimes W_{3})^{*}$$
is defined by
$$(\Psi_{P(z_{1}, z_{2})}^{(1)}(\nu))
(w_{1}\otimes w_{2}\otimes w_{3})=\langle \nu,
w_{1}\boxtimes_{P(z_{1})}
(w_{2}\boxtimes_{P(z_{2})} w_{3})\rangle$$
for $\nu\in W_{1}\hboxtr_{P(z_{1})}(W_{1}\boxtimes W_{3})$,
$w_{1}\in W_{1}$, $w_{2}\in W_{2}$ and $w_{3}\in W_{3}$.
In addition, by the properties of (\ref{prod1}) and the linear
independence
of $\left(\log\left(1-\frac{z_{2}}{z_{1}}\right)\right)^{k}$,
$k=1, \dots, K$, $g_{K}(w_{0}, w_{1}, w_{2}, w_{3};
z_{1}, z_{2})$ also satisfy the
$L(-1)$-derivative properties
\begin{eqnarray}
\frac{\partial}{\partial z_{1}}g_{K}(w_{0}, w_{1}, w_{2}, w_{3};
z_{1}, z_{2})&=&g_{K}(w_{0}, L(-1)w_{1}, w_{2}, w_{3};
z_{1}, z_{2}),\label{gkl-1-1}\\
\frac{\partial}{\partial z_{2}}g_{K}(w_{0}, w_{1}, w_{2}, w_{3};
z_{1}, z_{2})&=&g_{K}(w_{0}, w_{1}, L(-1)w_{2}, w_{3};
z_{1}, z_{2}).\label{gkl-1-2}
\end{eqnarray}
Thus we see that $g_{K}(w_{0}, w_{1}, w_{2}, w_{3};
z_{1}, z_{2})$
can also be written in the form of (\ref{prod1}) with possibly
different
$\mathcal{Y}_{1}$ and $\mathcal{Y}_{2}$
in the region $|z_{1}|>|z_{2}|>0$.

By the
properties of (\ref{prod1}) and the linear independence
of $\left(\log\left(1-\frac{z_{2}}{z_{1}}\right)\right)^{k}$,
$k=1, \dots, K$, $g_{K-1}(w_{0}, w_{1}, w_{2}, w_{3};
z_{1}, z_{2})$ also satisfies all the properties satisfied by
(\ref{prod1}) except for the $L(-1)$-derivative property.
Thus for any fixed complex numbers $z_{1}$ and $z_{2}$ satisfying
$|z_{1}|>|z_{2}|>0$, using Conditions 1--3, Theorem \ref{fusion},
Theorem 14.10 in \cite{H1},
the construction of the associativity isomorphisms
in \cite{H1} and the same argument as the one above for
$g_{K}(w_{0}, w_{1}, w_{2}, w_{3};
z_{1}, z_{2})$, we see that
$g_{K-1}(w_{0}, w_{1}, w_{2}, w_{3};
z_{1}, z_{2})$ can also be written in the form of (\ref{prod1})
with possibly different $\mathcal{Y}_{1}$ and $\mathcal{Y}_{2}$.
Using  Conditions 1--2 and Theorem \ref{fusion},
we can always find finite basis of spaces of intertwining operators
among irreducible $V$-modules such that the products of these basis
intertwining operators
give us a finite basis of the space of linear functional on
$W_{0}\otimes W_{1}\otimes W_{3}\otimes W_{4}$ given by
products of intertwining operators of the form
(\ref{prod1}) (evaluated at the given numbers $z_{1}$ and $z_{2}$).
We denote the value of these basis elements of linear functionals
at $w_{0}\otimes w_{1}\otimes w_{2}\otimes w_{3}$ for $w_{i}\in
W_{i}$,
$i=0, 1, 2, 3$, by $e_{l}(w_{0}, w_{1}, w_{2}, w_{3}; z_{1}, z_{2})$,
$l=1, \dots, p$. Then we have
\begin{equation}\label{g-k-1}
g_{K-1}(w_{0}, w_{1}, w_{2}, w_{3};
z_{1}, z_{2})=\sum_{l=1}^{p}c_{l}(z_{1}, z_{2})
e_{l}(w_{0}, w_{1}, w_{2}, w_{3}; z_{1}, z_{2}),
\end{equation}
where we have written down the dependence of the coefficients
$c_{l}(z_{1}, z_{2})$, $l=1, \dots, p$,
on $z_{1}$ and $z_{2}$ explicitly. Since
$g_{K-1}(w_{0}, w_{1}, w_{2}, w_{3};
z_{1}, z_{2})$ and $e_{l}(w_{0}, w_{1}, w_{2}, w_{3}; z_{1}, z_{2})$,
$l=1, \dots, p$, are all analytic in $z_{1}$ and $z_{2}$, we see that
$c_{l}(z_{1}, z_{2})$, $l=1, \dots, p$, are also analytic
in $z_{1}$ and $z_{2}$. The functions
$e_{l}(w_{0}, w_{1}, w_{2}, w_{3}; z_{1}, z_{2})$ for
$l=1, \dots, p$ also
satisfy the $L(-1)$-derivative properties
\begin{eqnarray}
\frac{\partial}{\partial z_{1}}e_{l}(w_{0}, w_{1}, w_{2}, w_{3};
z_{1}, z_{2})&=&e_{l}(w_{0}, L(-1)w_{1}, w_{2}, w_{3};
z_{1}, z_{2}),\label{ell-1-1}\\
\frac{\partial}{\partial z_{2}}e_{l}(w_{0}, w_{1}, w_{2}, w_{3};
z_{1}, z_{2})&=&e_{l}(w_{0}, w_{1}, L(-1)w_{2}, w_{3};
z_{1}, z_{2}).\label{ell-1-2}
\end{eqnarray}

{}From the $L(-1)$-derivative property for intertwining operators
and the linear independence
of $\left(\log\left(1-\frac{z_{2}}{z_{1}}\right)\right)^{k}$,
$k=1, \dots, K$, we obtain
\begin{eqnarray}
\lefteqn{\frac{\partial}{\partial z_{1}}g_{K-1}(w_{0}, w_{1}, w_{2},
w_{3};
z_{1}, z_{2})}\nn
&&+K\frac{\partial}{\partial z_{1}}\log
\left(1-\frac{z_{2}}{z_{1}}\right)
g_{K}(w_{0}, w_{1}, w_{2}, w_{3};
z_{1}, z_{2})\nn
&&\quad =g_{K-1}(w_{0}, L(-1)w_{1}, w_{2}, w_{3};
z_{1}, z_{2}),\label{gk-1l-1-1}\\
\lefteqn{\frac{\partial}{\partial z_{2}}g_{K-1}(w_{0}, w_{1},
w_{2}, w_{3};
z_{1}, z_{2})}\nn
&&+K\frac{\partial}{\partial z_{2}}\log
\left(1-\frac{z_{2}}{z_{1}}\right)
g_{K}(w_{0}, w_{1}, w_{2}, w_{3};
z_{1}, z_{2})\nn
&&\quad =g_{K-1}(w_{0}, w_{1}, L(-1)w_{2}, w_{3};
z_{1}, z_{2}).\label{gk-1l-1-2}
\end{eqnarray}
{}From (\ref{gkl-1-1}), (\ref{gkl-1-2}),
(\ref{gk-1l-1-1}) and (\ref{gk-1l-1-2}), we obtain
\begin{eqnarray}
\lefteqn{\frac{\partial^{2}}{\partial z^{2}_{1}}
g_{K-1}(w_{0}, w_{1}, w_{2}, w_{3};
z_{1}, z_{2})-2\frac{\partial}{\partial z_{1}}
g_{K-1}(w_{0}, L(-1)w_{1}, w_{2}, w_{3};
z_{1}, z_{2})}\nn
&&\quad  \quad\quad  \quad\quad
\quad\quad+g_{K-1}(w_{0}, L(-1)^{2}w_{1}, w_{2}, w_{3};
z_{1}, z_{2})\nn
&&\quad  \quad\quad \quad\quad\quad\quad
+\left(\frac{\partial}{\partial z_{1}}\log
\left(1-\frac{z_{2}}{z_{1}}\right)\right)^{-1}
\frac{\partial^{2}}{\partial z_{1}^{2}}\log
\left(1-\frac{z_{2}}{z_{1}}\right)\cdot\nn
&&\quad  \quad \quad  \quad\quad \quad\quad\quad\quad
\cdot g_{K-1}(w_{0}, L(-1)w_{1}, w_{2}, w_{3};
z_{1}, z_{2})\nn
&&\quad  \quad\quad \quad\quad\quad\quad
-\left(\frac{\partial}{\partial z_{1}}\log
\left(1-\frac{z_{2}}{z_{1}}\right)\right)^{-1}
\frac{\partial^{2}}{\partial z_{1}^{2}}\log
\left(1-\frac{z_{2}}{z_{1}}\right)
\cdot\nn
&&\quad  \quad \quad  \quad\quad \quad\quad\quad\quad
\cdot \frac{\partial}{\partial z_{1}}g_{K-1}(w_{0}, w_{1}, w_{2},
w_{3};
z_{1}, z_{2})\nn
&&\quad  \quad\quad \quad\quad\quad\quad =0,\label{g-k-1-eqn1}\\
\lefteqn{\frac{\partial^{2}}{\partial z^{2}_{2}}
g_{K-1}(w_{0}, w_{1}, w_{2}, w_{3};
z_{1}, z_{2})-2\frac{\partial}{\partial z_{2}}
g_{K-1}(w_{0}, w_{1}, L(-1)w_{2}, w_{3};
z_{1}, z_{2})}\nn
&& \quad  \quad\quad \quad\quad
\quad\quad+g_{N-1}(w_{0}, w_{1}, L(-1)^{2}w_{2}, w_{3};
z_{1}, z_{2})\nn
&&\quad  \quad\quad \quad\quad\quad\quad
+\left(\frac{\partial}{\partial z_{2}}\log
\left(1-\frac{z_{2}}{z_{1}}\right)\right)^{-1}
\frac{\partial^{2}}{\partial z_{2}^{2}}\log
\left(1-\frac{z_{2}}{z_{1}}\right)\cdot\nn
&&\quad  \quad \quad  \quad\quad \quad\quad\quad\quad
\cdot
g_{K-1}(w_{0}, w_{1}, L(-1)w_{2}, w_{3};
z_{1}, z_{2})\nn
&&\quad  \quad\quad \quad\quad\quad\quad
-\left(\frac{\partial}{\partial z_{1}}\log
\left(1-\frac{z_{2}}{z_{1}}\right)\right)^{-1}
\frac{\partial^{2}}{\partial z_{1}^{2}}\log
\left(1-\frac{z_{2}}{z_{1}}\right)
\cdot\nn
&&\quad  \quad\quad  \quad\quad \quad\quad\quad\quad
\cdot \frac{\partial}{\partial z_{2}}g_{K-1}(w_{0}, w_{1}, w_{2},
w_{3};
z_{1}, z_{2})\nn
&&\quad  \quad\quad \quad\quad\quad\quad=0.\label{g-k-1-eqn2}
\end{eqnarray}

Substituting (\ref{g-k-1}) into (\ref{g-k-1-eqn1}) and
(\ref{g-k-1-eqn2}) and then using the $L(-1)$-derivative
properties (\ref{ell-1-1}) and (\ref{ell-1-2})
and the linear independence of
$e_{l}(w_{0}, w_{1}, w_{2}, w_{3}; z_{1}, z_{2})$ for
$l=1, \dots, p$, we obtain
\begin{eqnarray*}
\frac{\partial^{2}}{\partial z^{2}_{1}}c_{l}(z_{1}, z_{2})
&=&\left(\frac{\partial}{\partial z_{1}}\log
\left(1-\frac{z_{2}}{z_{1}}\right)\right)^{-1}
\frac{\partial^{2}}{\partial z_{1}^{2}}\log
\left(1-\frac{z_{2}}{z_{1}}\right)
\frac{\partial}{\partial z_{1}}c_{l}(z_{1}, z_{2}),\\
\frac{\partial^{2}}{\partial z^{2}_{2}}c_{l}(z_{1}, z_{2})
&=&\left(\frac{\partial}{\partial z_{2}}\log
\left(1-\frac{z_{2}}{z_{1}}\right)\right)^{-1}
\frac{\partial^{2}}{\partial z_{2}^{2}}\log
\left(1-\frac{z_{2}}{z_{1}}\right)
\frac{\partial}{\partial z_{2}}c_{l}(z_{1}, z_{2})
\end{eqnarray*}
for $l=1, \dots, p$.
The general solution of this system of equations
is
$$c_{l}(z_{1},
z_{2})=\lambda^{(l)}_{1}\log\left(1-\frac{z_{2}}{z_{1}}\right)
+\lambda^{(l)}_{2},\;\;\;l=1, \dots, p,$$
where $\lambda^{(l)}_{1}$ and $\lambda^{(l)}_{2}$,
$l=1, \dots, p$, are constants. Since
$g_{K-1}(w_{0}, w_{1}, w_{2}, w_{3};
z_{1}, z_{2})$ cannot contain terms proportional to
$\log\left(1-\frac{z_{1}}{z_{2}}\right)$, we have
$\lambda^{(l)}_{1}=0$ for $l=1, \dots, p$. So
$c_{l}(z_{1}, z_{2})$ for $l=1, \dots, p$ are constants. Thus
we have
$$\frac{\partial}{\partial z_{1}}g_{K-1}(w_{0}, w_{1}, w_{2}, w_{3};
z_{1}, z_{2})=g_{K-1}(w_{0}, L(-1)w_{1}, w_{2}, w_{3};
z_{1}, z_{2})$$
which contradicts to (\ref{gk-1l-1-1}). So $K=0$.

We now prove that there exists $N$ such that (\ref{si}) holds.

For any
$w_{0}\in W_{0}$, $w_{3}\in W_{3}$ and $z_{1}, z_{2}\in \mathbb{C}$
satisfying $|z_{1}|>|z_{2}|>|z_{1}-z_{2}|>0$, let
$\mu_{w_{0}, w_{3}}^{z_{1}, z_{2}}\in (W_{1}\otimes W_{2})^{*}$
be defined by
$$\mu_{w_{0}, w_{3}}^{z_{1}, z_{2}}(w_{1}\otimes w_{2})
=\langle w_{0}, \mathcal{Y}_{1}(w_{1}, x_{1})
\mathcal{Y}_{2}(w_{2}, x_{2})
w_{3}\rangle|_{x^{r}_{1}=e^{r\log z_{1}}, x_{2}^{r}=e^{r\log z_{2}},
r\in
\mathbb{R}}$$
for $w_{1}\in W_{1}$ and $w_{2}\in W_{2}$.
A straightforward calculation shows that
$\mu_{w_{0}, w_{3}}^{z_{1}, z_{2}}$ satisfies the
$P(z_{1}-z_{2})$-compatibility condition (cf. \cite{H1}).

We have proved that
\begin{equation}\label{lb1}
\langle w_{0}, \mathcal{Y}_{1}(w_{1}, x_{1})
\mathcal{Y}_{2}(w_{2}, x_{2})
w_{3}\rangle
=\sum_{i=1}^{j}z_{2}^{r_{i}}(z_{1}-z_{2})^{s_{i}}
f_{i}\left(\frac{z_{1}-z_{2}}{z_{2}}\right).
\end{equation}
Expanding $f_{i}$, $i=1, \dots, j$, we can write the right-hand side
of
(\ref{lb1}) as
$$\sum_{i=1}^{j}\sum_{m\in \mathbb{N}}C_{i, m}(w_{0}, w_{1}, w_{2},
w_{3})
z_{2}^{r_{i}-m}(z_{1}-z_{2})^{s_{i}+m}.$$
For $i=1, \dots, j$, $m\in \mathbb{N}$,
$w_{0}\in W_{0}$ and $w_{3}\in W_{3}$, let
$\beta_{i, m}(w_{0}, w_{3})\in (W_{1}\otimes W_{3})^{*}$ be defined
by
$$(\beta_{i, m}(w_{0}, w_{3}))(w_{1}\otimes w_{2})
=C_{i, m}(w_{0}, w_{1}, w_{2}, w_{3})$$
for $w_{1}\in W_{1}$ and $w_{2}\in W_{2}$.
Then we have
$$\mu_{w_{0}, w_{3}}^{z_{1}, z_{2}}=\sum_{i=1}^{j}\sum_{m\in
\mathbb{N}}
\beta_{i, m}(w_{0}, w_{3})e^{(r_{i}-m)\log z_{2}}
e^{(s_{i}+m)\log (z_{1}-z_{2})}.$$
Since $\mu_{w_{0}, w_{3}}^{z_{1}, z_{2}}$ satisfies the
$P(z_{1}-z_{2})$-compatibility condition,
$\beta_{i, m}(w_{0}, w_{3})$ for $i=1, \dots, j$ and $m\in \mathbb{N}$
also satisfy the
$P(z_{1}-z_{2})$-compatibility condition.

A straightforward calculation shows that
$\beta_{i, m}(w_{0}, w_{3})$ for $i=1, \dots, j$ and $m\in \mathbb{N}$
are eigenvectors of $L'_{P(z_{1}-z_{2}}(0)$ with eigenvalues
$\wt w_{1}+\wt w_{2}+s_{i}+m$ (cf. \cite{H1}).
For fixed $i=1, \dots, j$ and $m\in \mathbb{N}$,
we consider the $\C$-graded generalized $V$-module generated
by $\beta_{i, m}(w_{0}, w_{3})$. By Condition 1,
this $\C$-graded generalized $V$-module is a direct sum of irreducible
$V$-modules. The generator $\beta_{i, m}(w_{0}, w_{3})$
must belong to a finite sum of these irreducible
$V$-modules and thus this $\C$-graded generalized $V$-module
is in fact a module. So $\beta_{i, m}(w_{0}, w_{3})
\in W_{1}\hboxtr_{P(z_{1}-z_{2})}W_{2}$ for $i=1, \dots, j$
and $m\in \mathbb{N}$. Let $N$ be an integer such that
$(W_{1}\hboxtr_{P(z_{1}-z_{2})}W_{2})_{(r)}=0$ when $r\le N$.
Since the weights of $\beta_{i, m}(w_{0}, w_{3})$ are
$\wt w_{1}+\wt w_{2}+s_{i}+m$ for $i=1, \dots, j$ and $m\in
\mathbb{N}$,
we must have $\wt w_{1}+\wt w_{2}+s_{i}>N$.

The absolute convergence of (\ref{p-prod}) is proved using induction
together with
the same argument as in the proof of the absolute
convergence of (\ref{prod1}) except that instead of
the systems (\ref{sys2}) and Theorem \ref{regular}, the systems in
Theorem \ref{sys2} and Theorem \ref{sys2-regular} are used.
The analytic extension is obtained immediately 
using the system  (\ref{sys2}).
\epfv

\begin{rema}
{\rm The proof of Theorem \ref{conv-ext} is in fact the main work
in this paper.
It uses the results obtained in Sections 1 and 2.  
Note 
that the proof of the convergence and the proof that the matrix 
elements of a 
product of intertwining operators 
can be analytically extended to a multivalued function 
of a certain form uses only the regularity 
of a singular point of ordinary differential equations
induced from the system (\ref{eqn1}) 
and (\ref{eqn2}). There is no need to use 
the regularity of the sets of singular points of the 
two-variables system (\ref{eqn1}) 
and (\ref{eqn2}) as proved and discussed in Remark \ref{regular-2-var}.
The remaining part of
the proof in this section
is to show that there are no logarithm terms.
This remaining part of the proof  uses
the version of the associativity theorem proved in \cite{H1}
stated in terms
of the $P(z_{1}, z_{2})$-compatibility condition, $P(z_{1},
z_{2})$-local grading-restriction condition and the $P(z_{2})$-local
grading-restriction condition. Note what we want to prove here is the
convergence and extension property. So the other more widely known
version of the associativity
theorem in \cite{H1} stated in terms of convergence and extension
property
cannot be used here. }
\end{rema}

The following  result is the main application
we are interested in the present paper:

\begin{thm}\label{application1}
Let $V$ be a vertex operator algebra satisfying the
three conditions in Theorem \ref{conv-ext}.
Then the direct sum of all (inequivalent) irreducible
$V$-modules has a natural structure of intertwining operator algebra
and
the category of $V$-modules has a natural structure of vertex tensor
category. In particular, the category of $V$-modules has a natural
structure
of braided tensor category.
\end{thm}
\pf
Combining Theorems \ref{conv-ext}, Theorem \ref{fusion} and
the results  in \cite{HL1}--\cite{HL5}, \cite{H1},
\cite{H3}--\cite{H4} and \cite{H6},
we obtain Theorem \ref{application1} immediately.
\epfv

\begin{rema}
{\rm Note that by the results obtained in \cite{Z},
\cite{FZ}, \cite{H1}, \cite{DLM1}--\cite{DLM2}, \cite{L}, \cite{GN},
\cite{B}, \cite{ABD} and \cite{HKL}, the conclusion of the theorem
above
holds if various other  useful conditions hold. Here we discuss one
example.
Let $V$ be a
vertex operator algebra. Let $C_{2}(V)$ be the subspace of
$V$ spanned by elements of the form $u_{-2}v$ for $u, v\in V$.
Recall that $V$ is said to be {\it $C_{2}$-cofinite} or
satisfy the {\it $C_{2}$-cofiniteness condition}
if $\dim
V/C_{2}(V)<\infty$. We assume that
$V_{(n)}=0$ for $n<0$ and $V_{(0)}=\mathbb{C}\mathbf{1}$,
every
$\mathbb{N}$-gradable weak $V$-module
is completely reducible and $V$ satisfies the $C_{2}$-cofiniteness
condition.
By the results of Abe, Buhl and Dong in \cite{ABD},
Condition 1 is a consequence of these three conditions and
Condition 3 is a consequence of the
$C_{2}$-cofiniteness condition for $V$. By a result of
Anderson-Moore
\cite{AM} and Dong-Li-Mason \cite{DLM3}, in this case,
every irreducible $\C$-graded $V$-module is
in fact $\mathbb{Q}$-graded. Thus in this case, the conclusion of
Theorem \ref{application1} holds.}
\end{rema}

In the example discussed in the
remark above, we still need the condition that
every $\mathbb{N}$-gradable weak $V$-module is completely
reducible.  Here we have a much stronger result:

\begin{thm}\label{application2}
Let $V$ be a vertex operator algebra satisfying the following
conditions:

\begin{enumerate}

\item For $n<0$, $V_{(n)}=0$  and $V_{(0)}=\mathbb{C}\mathbf{1}$.

\item Every $\C$-graded $V$-module is completely reducible.

\item $V$ is $C_{2}$-cofinite.

\end{enumerate}
Then every $\C$-graded $V$-module is $C_{1}$-cofinite and
every $\C$-graded generalized $V$-module is a direct sum of
irreducible $\C$-graded $V$-modules. In particular, if $V$ satisfies
these three conditions  and every  (irreducible) $\C$-graded
$V$-module is $\mathbb{R}$-graded, then the conclusions of
Theorems \ref{conv-ext} and  \ref{application1} hold.
\end{thm}
\pf
Since $V$ satisfies Condition 3, Zhu's algebra $A(V)$ is
finite-dimensional (see \cite{DLM3}). Thus there are only finitely
many irreducible $A(V)$-modules. Since there is a bijection between
the set of equivalence classes of
irreducible $A(V)$-modules and the set of equivalence classes of
irreducible $\C$-graded $V$-modules (see \cite{Z}),
we see that there are only
finitely many inequivalent irreducible $\C$-graded $V$-modules.

Since $V$ satisfies Conditions 1 and 3,
by Corollary 3.18 in \cite{KarL},
every finitely-generated lower-truncated $\C$-graded
generalized $V$-module is a
$\C$-graded $V$-module.

The next part of the proof uses the results of Abe, Buhl and Dong
\cite{ABD} and, for a large part, is similar to the proof of some
results in \cite{ABD}.

By Condition 2, to prove that every $\C$-graded
$V$-module is $C_{1}$-cofinite, we need
only consider irreducible $\C$-graded $V$-modules. By a result in
\cite{ABD},
any irreducible $\C$-graded $V$-module is $C_{2}$-cofinite.
Since $C_{2}$-cofinite $\C$-graded
$V$-module is $C_{1}$-cofinite when $V$ satisfies Condition 1,
every $\C$-graded $V$-module is $C_{1}$-cofinite.

In \cite{ABD}, it was proved that if $V$ satisfies Conditions
1 and 3,
then any weak $V$-module $W$ has a nonzero lowest-weight vector,
that is, a vector $w\in W$ such that $v_{n}w=0$ when $\wt v_{n}<0$.
In particular, any $\C$-graded generalized $V$-module has a nonzero
lowest-weight vector. Clearly the $\C$-graded generalized
$V$-submodule
generated by $w$ is a finitely-generated
lower-truncated $\C$-graded generalized $V$-module. Since we have
proved that
such $\C$-graded generalized $V$-modules are $\C$-graded $V$-modules,
the $\C$-graded generalized $V$-submodule
generated by $w$ is a $\C$-graded $V$-module.
So we see that any nonzero $\C$-graded generalized $V$-module
contains a non-zero $V$-submodule. By Condition 2,
we see that any nonzero $\C$-graded generalized $V$-module in fact
contains
a nonzero  irreducible $\C$-graded $V$-submodule.

Now a standard argument allows us to show that
every $\C$-graded generalized $V$-module is a direct sum of
irreducible $\C$-graded $V$-modules. In fact,
given any $\C$-graded generalized $V$-module $W$, let $\tilde{W}$ be
the
sum of all irreducible $\C$-graded $V$-modules contained in $W$. Since
we have
proved that there are only finitely many irreducible $\C$-graded
$V$-modules,
$\tilde{W}$ is a lower-truncated $\C$-graded generalized $V$-module.
If $\tilde{W}\ne W$, $W/\tilde{W}$ is a nonzero $\C$-graded
generalized $V$-module
and thus contains a nonzero irreducible $\C$-graded
$V$-submodule $W_{0}/\tilde{W}$.
Since both $W_{0}/\tilde{W}$ and
$\tilde{W}$ are lower-truncated $\C$-graded generalized $V$-modules,
$W_{0}$ is also such a $\C$-graded generalized $V$-module.
It is clear  that Condition 2 and the fact that
every finitely-generated lower-truncated $\C$-graded
generalized $V$-module is a $\C$-graded
$V$-module are equivalent to the fact that every
lower-truncated $\C$-graded generalized $V$-module is a direct
sum of irreducible $\C$-graded $V$-modules. So $W_{0}$ is a direct sum
of
irreducible $\C$-graded $V$-modules. Since $W_{0}$ contains
$\tilde{W}$ and
$W_{0}/\tilde{W}$ is not $0$, $W_{0}$ as a direct
sum of irreducible $\C$-graded $V$-module
contains more irreducible $\C$-graded
$V$-submodules of $W$ than $\tilde{W}$.
Contradiction.

The last conclusion follows from the first two conclusions we
have just proved and
Theorems \ref{conv-ext} and \ref{application1}.
\epfv

\noindent {\small \sc Department of Mathematics, Rutgers University,
110 Frelinghuysen Rd., Piscataway, NJ 08854-8019 }

\noindent {\em E-mail address}: yzhuang@math.rutgers.edu

\end{document}